\documentclass{article}

\usepackage{arxiv}
\RequirePackage{amsthm,amsmath,amsfonts,amssymb}

\usepackage[utf8]{inputenc} 
\usepackage[T1]{fontenc}    
\usepackage{hyperref}       
\usepackage{url}            
\usepackage{booktabs}       
\usepackage{amsfonts}       
\usepackage{nicefrac}       
\usepackage{microtype}      
\usepackage{cleveref}       
\usepackage{lipsum}         
\usepackage{graphicx}
\usepackage{doi}

\usepackage{graphicx}
\usepackage[shortlabels]{enumitem}
\usepackage{float}
\usepackage{subcaption}
\usepackage[numbers]{natbib}

\theoremstyle{plain}

\newtheorem{theorem}{Theorem}[section]
\newtheorem{lemma}[theorem]{Lemma}
\newtheorem{corollary}{Corollary}[theorem]

\newtheorem{proposition}{Proposition}[theorem]

\theoremstyle{remark}

\newtheorem*{remark}{Remark}

\title{Non-stationary Lattice Anderson Model with Non-local Laplacian and Correlated White Noise}

\author{ 
   Xiaoyun Chen\\
   Department of Mathematics\\
   University of North Carolina at Charlotte\\
   9201 University City Blvd, Charlotte, NC,28223,USA\\
   \And 
   \href{https://orcid.org/0000-0002-1253-9076}{\includegraphics[scale=0.06]{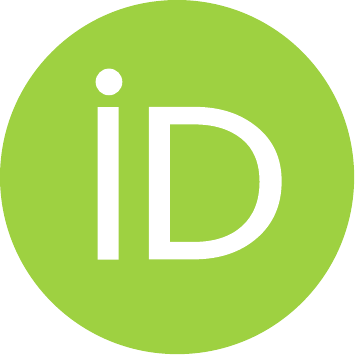}\hspace{1mm}Dan Han}\thanks{Correspondence Author} \\
    Department of Mathematics\\
    University of Louisville\\
    2301 S 3rd St, Louisville, KY 40292, USA\\
    	\texttt{dan.han@louisville.edu}\\
	\And
	Stanislav Molchanov \\
	Department of Mathematics\\
	University of North Carolina at Charlotte\\
	9201 University City Blvd, Charlotte, NC,28223,USA\\
	\texttt{smolchan@uncc.edu} \\
}

\renewcommand{\shorttitle}{\textit{arXiv} Template}

\hypersetup{
pdftitle={A template for the arxiv style},
pdfsubject={q-bio.NC, q-bio.QM},
pdfauthor={David S.~Hippocampus, Elias D.~Striatum},
pdfkeywords={First keyword, Second keyword, More},
}

\begin{document}
\maketitle

\begin{abstract}
We study the non-stationary Anderson parabolic problem on the lattice $Z^d$, i.e., the equation 
\begin{equation}\label{andersonmodel}
  \begin{aligned}
      \frac{\partial u}{\partial t} &=\varkappa \mathcal{A}u(t,x)+\xi_{t}(x)u(t,x)\\
       u(0,x) &\equiv 1, \, (t,x) \in [0,\infty)\times Z^d.
  \end{aligned}
\end{equation}
Here $\mathcal{A}$ is non-local Laplacian, $\xi_t (x), \ t \geq 0, \ x \in Z^d$ is the family of the correlated white noises and $\varkappa >0$ is the diffusion coefficient. The changes of $\varkappa$ (large versus small) are responsible for the qualitative phase transition in the model. At the first step the analysis of the model is reduced to the solution of the stochastic differential equation(SDE) (in the standard It\^{o}'s form) on the weighted Hilbert space $l^2(Z^d,\mu)$ with appropriate measure $\mu$. The equations of first two moments of the solution $u(t,x)$ are derived and studied using the spectral analysis of the corresponding  Schr\"{o}dinger operators with special class of the positive definite potentials. The analysis reveals several bifurcations depending on the properties of the kernel of $\mathcal{A}$ and the correlation function in the potential.\end{abstract}


\section{Introduction}
In this paper we'll develop several new results inspired by the memoir written by R. Carmona and S. Molchanov \cite{carmona1994parabolic}. Consider the non-stationary parabolic equation
\begin{equation}
  \begin{aligned}
      \frac{\partial u}{\partial t} &=\varkappa \mathcal{A}u(t,x)+\xi_{t}(x)u(t,x)\\
       u(0,x) &\equiv 1, \ (t,x) \in [0,\infty)\times Z^d.
  \end{aligned}
\end{equation}

Here the operator $\mathcal{A}$ is generator of the random walk $x(t)$ with continuous time. It has the form 
\begin{equation}
\mathcal{A}f(t,x)= \sum\limits_{z \neq 0, 
                  z \in Z^d }
                a(z)\big (f(t,x+z)-f(t,x)\big ).
\end{equation}
 We assume that $a(z)=a(-z)$(symmetry), $a(z)\geq 0$, $a(0)=0$, $a(z)>0$ if $|z|=1$(it gives non-periodicity),  $\sum\limits_{ z \in Z^d }a(z)=1$(normalization). The constant $\varkappa$ in (\ref{andersonmodel}) is the diffusion coefficient(diffusivity). The random walk has the following structure. It spends in each site $x \in Z^d$ the time $\tau_x$ which is exponential distributed with coefficient $\varkappa$, i.e., $P(\tau_x >t) =  e^{-\varkappa t}$ and at the moment $\tau_x+0$ it jumps from site $x$ to site $x+z$ with probability $a(z)$. We'll discuss the following three different cases. 
\begin{enumerate}[I)]
\item \textbf{Light Tails}
If $a(z)$ satisfies Cram\'{e}r's condition: 
\begin{equation}\label{light tails condition}
a(z) \leq ce^{-\eta|z|},\  c, \eta >0,
\end{equation}

then we'll say that the random walk has light tail. Then for the transition probability $p(t,0,z)$, one can use limit theorems with large deviations of Cram\'{e}r type. They give the complete control of $p(t,0,z) $ as $ t \to \infty$. There are many papers containing this information, among recent see \cite{Yarovaya}. This paper \cite{Yarovaya} contains important asymptotics of the Green Functions  $G_{\lambda}(0,z)= \int^{\infty}_{0} e^{\lambda t}p(t,0,z)dt$ where $\lambda>0$.

\item \textbf{Moderate Tails} 
\begin{equation}\label{moderate tails condition}
a(z) \sim  \frac{C(\dot{z})}{|z|^{d+ \alpha}},  \vert z \vert \to \infty, \  \dot{z} = \frac{z}{|z|} (direction\  of\  z)\in S ^{d-1},\  \alpha> 2,
\end{equation}
where $C$ is positive continuous function on $S ^{d-1}$, 
  $C(\dot{z})=C(-\dot{z})$. 
In particular, if $\sum|z|^2 a(z)< \infty$, the process $x(t)$ satisfies the central limit theorem, i.e., it has asymptotically Gaussian distribution with global expansion of the remainder term under some technical conditions, see \cite{molchanov2019population}.

\item \textbf{Heavy Tails}
\begin{equation}\label{heavy tails condition}
a(z)\sim \frac{C_{0}(\dot{z})}{|z|  ^{d+\alpha}},\,\,\vert z \vert \to \infty, \dot{z} = \frac{z}{|z|} (direction\  of\  z) \in S ^{d-1},\,\,0<\alpha<2
\end{equation}
where $C_{0}$ is positive continuous function on $S ^{d-1}$, 
  $C_{0}(\dot{z})=C_{0}(-\dot{z})$. 
  Under these conditions, the random variable $x(t)$ asymptotically belongs to the domain of attraction of the class of the $d$-dimensional symmetric stable laws with parameter $\alpha$. Under stronger conditions (the asymptotic expansion: $a(z)=\frac{C_{0}(\dot{z})}{\mid z\mid  ^{d+\alpha}}+\frac{C_{1}(\dot{z})}{\mid z\mid  ^{d+\alpha+1}}+\dots$  ) see \cite{getan2017intermittency}, one can prove global local limit theorem which gives the asymptotics of the transition probabilities $p(t,0,z)= p(t,x,x+z)= P \lbrace {x(t)= x+z \mid x(0)=x}\rbrace$ for $t \to \infty$ acting uniformly on $z \in Z^d$, see  \cite{getan2017intermittency}.

\end{enumerate}
 
We will use two different points of view on the operator $\mathcal{A}$. In the study of the stochastic differential equations in section 2,3,4 and 5, the operator $\mathcal{A}$ will be considered as the bounded operator in the weighted Hilbert space $l^2(Z^d,\mu)$ with the dot product 
$$
(f,g)_{\mu}= \sum\limits _{x \in Z^d}f(x)\ \bar{g}(x) \mu(x)
$$
where $\mu(x)>0$, $\sum\limits _{x \in  Z^d} \mu(x) < \infty$ with some additional technical conditions. 

In the spectral analysis of the Schr\"{o}dinger operators associated to the moments of the field $u(t,x)$, we use the standard space $l^2(Z^d)$ with the dot product 
$
(f,g)_{\mu}= \sum\limits _{x \in Z^d}f(x)\ \bar{g}(x) 
$. Operator $\mathcal{A}$ now is the bounded self-adjoint operator in $l^2(Z^d)$ and its perturbation, say, $\mathcal{H}_2=2\varkappa\mathcal{A}+B(x)$ by the positive definite function $B$ (the correlation of the potential $\xi_t(x)$) is the central object of the analysis in section 6, 7 and 8. 
 
Now we'll describe the potential $\xi_{t}(x)$. Formally, $\xi_{t}(x) = \dot{W}(t,x)$ where $\{W(t,x), t\geq 0, x\in Z^d\}$  is the family of Wiener processes at time $t\geq 0$ in location $x\in Z^d$.  Those processes have independent increments in time $t$ and represent the stationary Gaussian field on $Z^d$.   We will use notation $<\cdot>$ as the expectation over the law of the field $W(\cdot,\cdot)$ or white noise $\dot{W}(t,x), x\in Z^d$ and $E[\cdot]$ as the expectation over the law of random walk associated with the generator $\varkappa\mathcal{A}$. We will assume that the expectation
$<W(t,x)>=0$, $<W(t_{1},x) W(t_{2},y)>=(t_{1} \wedge t_{2})B(x-y)$ and $W(t,\cdot)\in l^2(Z^d,\mu)$, i.e. 
\begin{equation}
    \sum_{x \in Z^d}  W^2 (t,x) \mu(x) =  \Vert W(t, \cdot) \Vert_{\mu} ^2 < \infty P-\text{a.s.}
\end{equation}

 We'll construct $W(t,x)$ as the linear transformation of the field $\{w(t,\cdot)$, $t\geqslant 0\}$ where $w(t,x)$, $x \in Z^d$ are independent and identically distributed (i.i.d) standard Brownian motions at time $t\geq 0$ and location $x\in Z^d$. Define $W(t,x) = \sum\limits_{z \in Z^d} b(x-z) w(t,z)$, where $b(x-z)$ is the weight function for $w(t,x)$ and $\sum\limits_{z\in Z^d}|b(z)|<\infty$. Then  $< W(t_{1},x) W(t_{2},y) >=(t_{1} \wedge t_{2})B(x-y)$, where $t_1\wedge t_2=min\{t_1,t_2\}$, $B(x-y)=\sum\limits_{z\in Z^d}b(x-z)b(y-z)=\sum\limits_{z\in Z^d}b(x-y-z)b(-z)$. The proposed condition $\sum\limits_{z\in Z^d}|b(z)|<\infty$ leads to the convergence of $ \sum \limits_{z \in Z^d}  b^2(z)$, hence, the series $W(t,x)$ converges almost surely by Kolmogorov's three-series theorem.

For fixed $x \in Z^d$, the process $W(t,x)$ as the function of time $t$ is the Gaussian process with independent increments and $< W(t,x)W(s,x)> = (s\wedge t) B(0)$, where $B(0) = \sum\limits_{z \in Z^d} b^2(z)$, that is, $W(t,\cdot)=\sqrt{B(0)}w(t,\cdot)$ in law. 

However, we can also consider the realization of the field $W(t, \cdot) = \{W(t,x),  \ x \in Z^d\}$ as the elements of weighted Hilbert space $l^2(Z^d, \mu)$. Construction of the appropriate measure $\mu$ depending on the tails $a(z),  \ |z| \rightarrow \infty$  and the proof of continuity of $W(t, \cdot)$ in $l^2(Z^d,  \mu)$ are given in section 3. It also contains the estimation of $\| \mathcal{A}\|_{\mu}$ essential for the existence-uniqueness theorem for the parabolic Anderson model with  potential $\dot{W}(t,x)$ in section 4. 

Section 5 has many common points with the memoir \cite{carmona1994parabolic} and based on the classical It\^{o} calculus, additional results include the Kac-Feynman representation  of the solution $u(t,x)$ of the Anderson parabolic problem (\ref{andersonmodel}). Section 6 studies the moment equations for 
\begin{equation} \label{p-momnet definiton}
m_p(t, x_1, \cdots, x_p) = <u(t,x_1) \cdots u(t,x_p) >, \ (x_1, \cdots, x_p) \in Z^{dp}, \  p = 1,2,\cdots
\end{equation}
These equations are typical multiparticle Schr\"odinger equations.
\begin{equation}
\begin{aligned}
&\frac{\partial m_p}{\partial t} = \varkappa \left( \sum_{i=1}^p \mathcal{A}_{x_i}\right)m_p + V_p(x_1, \cdots, x_p)m_p\\
&m_p(0, x_1, \cdots, x_p)=1.
\end{aligned}
\end{equation}
Here $\mathcal{A}_{x_i}m_p= \sum\limits_{z \neq 0, 
                  z \in Z^d }
                a(z)\big (m_p(t,x_1,\cdots,x_i+z,\cdots,x_p)-m_p(t,x_1,\cdots,x_i,\cdots,x_p)\big )$, $i=1,2,\cdots,p$.  $
V_p(x_1, \cdots, x_p) =\sum\limits_{i <j}B(x_i  - x_j)
$ is $p-$particle potential with binary interaction $B(x_i - x_j) = <W(1, x_i) W(1, x_j)>$ . Derivation of moment equations is based on It\^{o} formula. 

The equation for the first moment $m_1(t,x)= <u(t,x)>$ is trivial:
\begin{equation}
\frac{\partial m_1}{\partial t}=\varkappa \mathcal{A}m_1, \ m_1(0,x) =1
\end{equation}
it gives $ m_1(t,x) \equiv 1 $.

The most important part of the work is devoted to the study of the moments $m_p(t,  \cdots),  p \geq 2$  and first of all to the second moments $m_2(t,  x_1-x_2) = <u(t,x_1)u(t,x_2)>$. Standard transition to the center of mass leads to equation 
\begin{equation} \label{2nd moments }
 \begin{aligned}
 &\frac{\partial m_2(t,z)}{\partial t} = 2 \varkappa \mathcal{A} m_2 + B(z)m_2 =\mathcal{H}_2 m_2 \\
 & m_2(0,z) \equiv 1.
  \end{aligned}
\end{equation}

Memoir \cite{carmona1994parabolic} contains the similar equation that has the form 
\begin{equation} \label{2nd moments with local operator }
\frac{\partial m_2(t,z)}{\partial t} = 2 \varkappa \Delta m_2 + \delta_0(x) m_2 = \tilde{\mathcal{H}_2} m_2. 
\end{equation}
However, the operator  $\Delta f(x)=\sum\limits_{x':|x-x'|=1}[f(x')-f(x)]$ is a local Laplacian instead of the nonlocal operator $\mathcal{A}$. This equation can be solved explicitly using Fourier transform. The central results of \cite{carmona1994parabolic} about equation (\ref{2nd moments with local operator }) shows : for $d =1,2$ operator in the right part of (\ref{2nd moments with local operator }) has positive eigenvlaue $\lambda_0(\varkappa) >0$ for any $\varkappa >0$. If $d \geq 3 $ there exists the bifurcation: $\lambda_0(\varkappa) >0$ exists for $\varkappa < \varkappa_{cr}$, if $\varkappa \geq \varkappa_{cr}$ the spectrum of $\tilde{\mathcal{H}}_2$ is pure a.c (there is an explicit formula for $\varkappa_{cr}$).  In other terms, the second moment of $u(t,x)$ is exponentially increasing for $d=1,2,  \ t \rightarrow \infty,  \varkappa >0$. If $d \geq 3$, the same is true only for small $ \varkappa(\varkappa < \varkappa_{cr})$.  If $d \geq 3 $ and $\varkappa  \geq \varkappa_{cr} $, then the second moment of $u(t,x)$ is bounded in time.

Positivity of ground state energy $\lambda_0$ is the manifestation of intermittency : high irregularity of the field $u(t,x)$; while  the absence of $\lambda_0 > 0$ demonstrates the regularity of the Anderson parabolic problem. 

There are many publications on the spectral theory of the lattice Schr\"odigner type operator mainly in the case of the local Laplacian, i.e,  $\Delta f(x)=\sum\limits_{x':|x-x'|=1}[f(x')-f(x)]$, although this area is not so well studied as the classical quantum mechanical Hamiltonian in $L^2(R^d)$. This paper is studying the nonlocal Schr\"odinger operator with correlation function (positive definite) $B$ as the potential. The estimation for the transition probability $p(t,x,y)$ for the random walk $x(t)$ associated with the nonlocal operator $2\varkappa\mathcal{A}$ are given in the section $7$. 
Section $8$  contains analysis of the additional spectral bifurcations for the operator $\mathcal{H}_2$. Such bifurcations are related to transition from the light tailed $a(\cdot)$ to the heavy tails, transition from the strong to weak correlation of $W(t,x), \  x \in Z^d$, transition from the case when $\sum\limits_{x \in Z^d} B(x) >0$ to $\sum\limits_{x \in Z^d} B(x) =0$ (note the case $\sum\limits_{x \in Z^d} B(x) <0$  is impossible) etc. 

Analysis of the Lyapunov exponents for the higher moment $m_p,  \  p \geq 3$ as well as $P$-a.s. Lyapunov exponents $\tilde{\gamma}(\varkappa)=\lim_{t \rightarrow \infty} \displaystyle\frac{\ln (u(t,x))}{t}$ will be the topic of the next publication.

\section{ Weighted Hilbert space $l^2(Z^d,\mu)$ and Properties of the non-local Laplacian $\mathcal{A}$ }

\subsection{Weighted Hilbert space $l^2(Z^d,\mu)$}
All the future analysis will be concentrated on the weighted Hilbert space $l^2(Z^d,\mu)=\{f(\cdot):\Vert f\Vert^2_{\mu}=\sum\limits_{x\in Z^d}|f(x)|^2\mu(x)\}$. 

\textbf{Regularity condition for weight $\mu$}:
 First of all, we'll introduce regularity condition on the weight $\mu(x)$, $x \in Z^d$. Let $\mu(x)=\mu(\vert x \vert)$, $\vert x \vert=\vert x \vert_{\infty}=\max\limits_{i}\vert x_i \vert$ and $\mu(n)$, $n=0,1,2,\cdots$ be the monotone concave down function of $n$. It means that
\begin{equation}\label{concave down of mu}
1 \leq \frac{\mu(0)}{\mu(1)} \leq \frac{\mu(1)}{\mu(2)} \leq \cdots \leq \frac{\mu(n)}{\mu(n+1)}.
\end{equation}
 Then 
\begin{equation}\label{mu_regular_condition0}
\frac{\mu(0)}{\mu(1)} \cdot \frac{\mu(1)}{\mu(2)} \cdots  \frac{\mu(n-1)}{\mu(n)} = \frac{\mu(0)}{\mu(n)} \leqslant  \frac{\mu(k)}{\mu(k+1)} \cdot \frac{\mu(k+1)}{\mu(k+2)} \cdots \frac{\mu(k+n-1)}{\mu(k+n)} = \frac{\mu(k)}{\mu(k+n)}.
\end{equation}

For fixed $n \geq 1$, $\sup\limits_{x,z:|x-z|=n}\displaystyle \frac{\mu(x)}{\mu(z)} =\frac{\mu(0)}{\mu(n)}$. Thus we get the following lemma

\begin{lemma}
For any $a \in Z^d$,
\begin{equation}\label{mu_regular_condition}
\sup_{x \in Z^d} \frac{\mu(x+a)}{\mu(x)} \leq h(|a|) = \frac{\mu(0)}{\mu(|a|)}.
\end{equation}
\end{lemma}

\subsection{Conditions for  $ \| \mathcal{A} \|_{\mu} < \infty$}

Our next goal is the estimation of the norm of $\mathcal{A}$ in $l^2(Z^d, \mu)$. In particular we'll describe all measure $\mu$ with restriction (\ref{mu_regular_condition}) such that $ \| \mathcal{A} \|_{\mu} < \infty$ in $l^2(Z^d, \mu)$. Instead of $\mathcal{A} f = \sum\limits_{z \neq 0, z\in Z^d} a(z)(f(t,x+z)-f(t,x))$, we consider $\bar{\mathcal{A}}f  = \sum\limits_{y \in Z^d} a(x-y)f(t,y) =\sum\limits_{z \in Z^d} a(z)f(t,x-z)$. One can note that $\| \mathcal{A} \|_{\mu} \leq  \| \bar{\mathcal{A}} \|_{\mu} +1$, thus if $\| \bar{\mathcal{A}} \|_{\mu}<\infty$, then $\| \mathcal{A} \|_{\mu}<\infty$.

\begin{lemma}\label{norm of A 1} Under the condition (\ref{mu_regular_condition0}),
 $$ \|\bar{\mathcal{A}} \|_\mu  \leq \sum_{z \neq 0, z\in Z^d} a(z) \sqrt{\frac{\mu(0)}{\mu(z)}}$$.
\end{lemma}  
\begin{proof}
Let $f(x) \in l^2(Z^d, \mu)$, then $ \|\bar{\mathcal{A}} f(t,x)\|_{\mu} \leq \sum\limits_{z \neq 0, z\in Z^d} a(z)\|f(t,x-z)\|_{\mu}$. 
\begin{align*}
& \|f(t,x-z)\|_{\mu}^2 \leq \sum_{z \in Z^d} f^2(x-z)\mu(x) \cdot \frac{\mu(x-z)}{\mu(x-z)} \leq \frac{\mu(0)}{\mu(z)} \|f\|_{\mu}^2\\
& \| \bar{\mathcal{A}}\|_\mu  \leq \sum_{z \neq 0, z\in Z^d} a(z) \sqrt{\frac{\mu(0)}{\mu(z)}}.
\end{align*}
\end{proof}
One can get the following better estimation than that in lemma \ref{norm of A 1}.
\begin{lemma}\label{norm of A 2}
Under the condition (\ref{mu_regular_condition0}),
$$ \| \bar{\mathcal{A}} \|_\mu  \leq  \sqrt{\sum_{z \neq 0, z\in Z^d}  \frac{ a^2(z)\mu(0)}{\mu(z)}}.$$
\end{lemma}

\begin{proof}

$$
\Vert \bar{\mathcal{A}}f(x)\Vert_{\mu} ^2 \leq \bigg(  \sum_{y\in Z^d} a(x-y) f(t,y)\bigg)^2 = \bigg(  \sum_{y} \frac{a(x-y)}{\sqrt{\mu(y)}}  f(t,y) \sqrt{\mu(y)} \bigg)^2  \leq \sum_{y} \frac{a^{2} (x-y)}{\mu(y)} \Vert f \Vert^{2}_{\mu}.
$$

Multiplying both parts by $\mu(x)$, we'll get
$$
\frac{\Vert \mathcal{A} f \Vert^{2}_{\mu}}{\Vert f \Vert^{2}_{\mu}}
 \leq \sum_{x,y} \frac{\mu(x)a^2(x-y)}{\mu(y)} = \sum_{x, \nu} \frac{\mu(y+z)a^2(z)}{\mu(y)} \leq \sum_{z} \frac{\mu(0)a^2(z)}{\mu(z)}. 
$$
It proves inequality (\ref{norm of A 2}). 

\end{proof}

\begin{theorem}\label{thm for boundness of A}

The following conditions are sufficient for the boundness of Laplacian $\mathcal{A}$ in $l^2(Z^d,\mu)$.

\begin{enumerate}
\item \textbf{Light Tail} 
If $a(z)$ satisfies Cram\'{e}r's condition:
\begin{equation}
a(z) \leq ce^{-\eta|z|},\  c, \eta >0,
\end{equation}
then
\begin{equation}
\| \bar{\mathcal{A}}\|_{\mu} \leq \sqrt{\sum_{z \in Z^d} c^2 e^{ - 2\eta |z|}\frac{\mu(0)}{\mu(z)}}.
\end{equation}
Thus if $ \mu(z) < e^{-\gamma |z|}$ and $ \gamma - 2 \eta<0$, then $\| \mathcal{A}\|_{\mu} < \infty$.

\item \textbf{Moderate Tail} or \textbf{heavy Tail} If $a(z)$ has moderate tail or heavy tail, then
\begin{equation}
a(z) \leq \frac{C(\dot{z})}{|z|^{d+ \alpha}},  \vert z \vert \to \infty, \  \dot{z} = \frac{z}{|z|} (direction\  of\  z)\in S ^{d-1} 
\end{equation}
where $C$ is positive continuous function on $S ^{d-1}$, $\alpha> 2$ for moderate tail, $\alpha\in (0,2)$ for heavy tail.
If $\mu(z)$ satisfies:
\begin{equation}
\mu(z) \sim \frac{C(\delta_2)}{1+|z|^{d+ \delta_2}},\,\,\,\delta_2>0.
\end{equation}
\begin{equation*}
\| \bar{\mathcal{A}}\|_{\mu} \leq \sqrt{\sum_{z \in Z^d} (a(z))^2 \frac{\mu(0)}{\mu(z)}}\sim \sqrt{\sum_{z \in Z^d} \frac{1+|z|^{d+\delta_2}}{|z|^{2d+2\alpha}}}=\sqrt{C_0+\sum_{z \in Z^d}\frac{1}{|z|^{d+2\alpha-\delta_2}}}
\end{equation*}
where $C_0$ is a finite constant. Then  $\| \mathcal{A}\|_{\mu} <\infty$ if $\delta_2<2\alpha$, i.e., for any sufficiently small  positive $\delta_2$.
\end{enumerate}
\end{theorem}

Such estimates are especially important in heavy tailed $a(z)$. We want to include in the theory transient random walks with the generator $\mathcal{A}$  in dimensions $d =1 , \ 2 $.

\section{Continuity of the process $W(t,x)$ in weighted Hilbert space $l^2(Z^d,\mu)$}
The Wiener process $W(t,x)$ is a weighted summation of  independent and identically distributed (i.i.d) standard Brownian motions $w(t,x)$, $x \in Z^d$ , $t\geq 0$: $W(t,x) = \sum\limits_{z \in Z^d} b(x-z) w(t,z)$. The series converges $P$-a.s. if and only if $\sum\limits_{z\in Z^d}b^2(z)<\infty$, that is $b(\cdot)\in l^2(Z^d)$. However, we will assume that the kernel $b(\cdot)$ is not only from $l^2(Z^d)$ but satisfies the stronger condition $b(\cdot) \in l^1(Z^d)$, $\sum\limits_{z\in Z^d}|b(z)| \leq \|b(\cdot)\|_1 < \infty$. Then
\begin{align*}
B(0) &= < W^2(1,0) > = \sum\limits_{z\in Z^d}b^2(z).\\
B(x) &= < W^2(1,x) > = \sum_{z \in Z^d} b(x-z)b(-z).
\end{align*}

i.e.,
$$
\sum_{x \in Z^d}|B(x)| \leq  \sum_{x, z \in Z^d} |b(x-z)| |b(z)| \leq \sum_{z \in Z^d} |b(z)|^2.
$$

Since $b(\cdot) \in l^1$, we have $b(x)=\displaystyle\frac{1}{(2\pi)^d}\int_{T^d}\hat{b}(k)e^{-ikx}dk$ where $\hat{b}(k)=\sum\limits_{x\in Z^d}b(x)e^{ikx}$. Then by Parseval's identity, $B(x)=\sum\limits_{z \in Z^d} b(x-z)b(-z)=\displaystyle\frac{1}{(2\pi)^d}\int_{T^d} e^{-ikx} |\hat{b}(k)|^2 dk$. Due to Bochner-Khinchin theorem, $B(x)$ is the positive definite function with continuous spectral density $\rho(k) =\displaystyle\frac{ |\hat{b}(k)|^2}{(2\pi)^d}$. Under the condition $b\in l^1(Z^d)$,
$$\sum\limits_{x \in Z^d}B(x)=\sum\limits_{x\in Z^d}\sum\limits_{z\in Z^d}b(x-z)b(-z)=(\sum\limits_{z\in Z^d}b(z))^2.$$

It means that $\sum\limits_{x\in Z^d} B(x)\geq 0$ and $\sum\limits_{x\in Z^d}B(x)=0$ if and only if $\sum\limits_{x\in Z^d}b(x)=0$. These facts will be crucial in the spectral analysis in the section 8. Note also that $\rho(0)=\frac{|\hat{b}(0)|^2}{(2\pi)^d}=0$ if and only if $\sum\limits_{x\in Z^d}b(x)=0$ or $\sum\limits_{x\in Z^d}B(x)=0$. This fact has the following probabilistic interpretation. 

Consider the Gaussian field $\xi(x)=W(1,x)$ and define $S_L=\sum\limits_{| x |_{\infty}\leq L}\xi(x)$ where $|x |_{\infty}=\max\limits_i |x_i|$ and $S_L^*=\frac{1}{(2L)^{d/2}}S_L$, then 
\begin{align*}
Var S_L^*&=\frac{1}{(2L)^{d}}<S_L^2>=\frac{1}{(2L)^d}\sum\limits_{x,y:\substack{| x|_{\infty}<L\\ | y|_{\infty}<L}} B(x-y)\\
&=\frac{1}{(2L)^d}\int_{T^d}|\sum\limits_{|x |_{\infty}<L}e^{ikx}|^2\rho(k)dk\\
&=\frac{1}{(2L)^d }\int_{T^d}\prod_{j=1}^d\frac{\sin^2((L+\frac{1}{2})k_j)}{\sin^2(\frac{k_j}{2})}\rho(k)dk.
\end{align*}
The Fejer kernel $\displaystyle\frac{1}{(2L)^d}\prod_{j=1}^d\frac{\sin(L+\frac{1}{2})k_j}{\sin(\frac{k_j}{2})}$ converges weakly in $C(T^d)$ to $(2\pi)^d\delta(k)$ if $L\rightarrow \infty$. The value of this density $\rho(k)$ at the point $k=0$ (long waves) is especially important. Namely, if $\rho(0) = 0$, then $\lim\limits_{L\rightarrow \infty}<\frac{S_L^2}{(2L)^d}>=0$ , but if $\rho(0) > 0$ then $\lim\limits_{L\rightarrow \infty}<\frac{S_L^2}{(2L)^d}>=(2\pi)^d\rho(0)=|\hat{b}(0)|^2$. For example, if $b(x)=\delta_0(x)$, $\hat{b}(k)=1$, then $\lim\limits_{L\rightarrow \infty}<\frac{S_L^2}{(2L)^d}>=1=|\hat{b}(0)|^2$.

\begin{theorem}
If $\sum\limits_{x \in Z^d} b^2(x) =B(0) < \infty$ and $\sum\limits_{z \in Z^d} \mu(x)< \infty$ , then Gaussian field $W(t,x) = \sum\limits_{z \in Z^d} b(x-z) w(t,z), t\geq 0, x \in Z^d$ as element of $l^2(Z^d, \mu)$ is a continuous function of time. 
\end{theorem}
\begin{proof}
Due to Kolmogorov's continuity criterion it is sufficiently to prove that for some $\alpha >0,  \  \delta>0,$
\begin{align*}
<\Vert \Delta W(\cdot, \cdot)\Vert_{\mu}^{\alpha}> \leq C |t_2 -t_1|^{1 + \delta}, \ \alpha >0,  \  \delta>0
\end{align*}
where $
\Vert \Delta W(\cdot, \cdot)\Vert_{\mu}^{\alpha} = (\sum_{x \in Z^d}\mu(x)((W(t_2, x)-W(t_1, x))^2 )^{\alpha/2}.
$
 
\begin{align*}
&<\Vert \Delta W\Vert_{\mu}^4>=<(\sum_{x \in Z^d}\mu(x)((W(t_2, x)-W(t_1, x))^2 )^2>\\
&=<\sum\limits_{x_1,x_2 \in Z^d}\mu(x_1)\mu(x_2)\Delta^2 W(\cdot,x_1)\Delta^2 W(\cdot,x_2)>
\end{align*}
where $\Delta W(\cdot,x_i)=\sum\limits_{z\in Z^d}b(x_i-z)(w(t_2,z)-w(t_1,z)),i=1,2$. 

The random variables $\Delta W(\cdot,x_1)$ and $\Delta W(\cdot,x_2)$ have the joint Gaussian distribution with parameters
\begin{align*}
<\Delta W(\cdot,x_1)>&=<\Delta W(\cdot,x_2)>=0,\\
<\Delta^2W(\cdot,x_1)>=<\Delta^2W(\cdot,x_2)>&=\sum\limits_{z\in Z^d}b^2(-z)|t_2-t_1|=B(0)|t_2-t_1|,\\
<\Delta W(\cdot,x_1)\Delta W(\cdot,x_2)>&=\sum\limits_{z\in Z^d}b(x_1-z)b(x_2-z)|t_2-t_1|=B(x_1-x_2)|t_2-t_1|.
\end{align*}
Then 
\begin{align*}
\displaystyle<e^{is_1\Delta W(\cdot,x_1)+is_2\Delta W(\cdot,x_2)}>=e^{-\displaystyle\frac{(s_1^2B(0)+2s_1s_2B(x_2-x_1)+s_2^2B(0))(t_2-t_1)}{2}}.
\end{align*}
Differentiation over $s_1$,$s_2$ for $s_1=s_2=0$ provides the fourth moments:
\begin{align*}
<\Delta^4 W(\cdot,x)>&=3B^2(0)|t_2-t_1|^2,\\
<\Delta^2 W(\cdot,x_1)\Delta^2 W(\cdot,x_2)>&=(B^2(0)+2B(0)B(x_2-x_1))|t_2-t_1|^2.
\end{align*}
The above result gives $<\Vert \Delta W(\cdot,\cdot)\Vert_{\mu}^4>\leq c|t_2-t_1|^2$, that is, the continuity of the functional valued process $W(t,x),x\in Z^d, t\in [0,T]$ in $l^2(Z^d,\mu)$.

One can study the higher moments of $\Delta W(t,\cdot)$ and find that $W(t,\cdot)$ belongs to any H\"older class $\mathbb{H}^{\beta}$ with $\beta< 1/2$ in $l^2(Z^d, \mu)$.
\end{proof}

\section{Existence and uniqueness of the solution $u(t,x)$}
We now study the problem of the existence and uniqueness of a solution of the model (\ref{whitenoisemodel}) 

\begin{equation}\label{whitenoisemodel}
  \begin{aligned}
      \frac{\partial u}{\partial t} &=\varkappa \mathcal{A}f(t,x)+\dot{W}(t,x)u(t,x)\\
       u(0,x) &\equiv 1, \ (t,x) \in [0,\infty)\times Z^d
  \end{aligned}
\end{equation}

where the potential $\xi_t(x)$ is $\dot{W}(t,x)$, and $W(t,x)=\sum\limits_{z\in Z^d}b(x-z)w(t,z)$ in which $w(t,z), z\in Z^d$ are independent standard Brownian motions. 
This equation can be considered either in usual It\^{o} form or in Stratonovich form which is popular in physical literatures. In physical literature, the white noise $\dot{w}_t$ is necessarily understood dynamically, as the sequence of the Gaussian process $\xi_{\varepsilon}(t)$, which correlation function $B_{\varepsilon}(t)$, converges to $\delta_0(t)$ if $\varepsilon \rightarrow 0$, see discussion in \cite{carmona1994parabolic}). We'll use the It\^{o} form of stochastic partial differential equations (SPDEs), i.e., we can understand  equation (\ref{whitenoisemodel}) as
\begin{equation}\label{usolution}
u(t,x)=1+\int_0^t\mathcal{A}u(s,x)ds+\int_0^tu(s,x)dW(s,x).
\end{equation}
\begin{theorem}
Assume that the white noise and the measure $\mu$ satisfy the conditions
$\sum\limits_{z\in Z^d}b^2(z)<\infty$, $\sum\limits_{z\in Z^d}\mu(x)<\infty$, and conditions mentioned in Theorem \ref{thm for boundness of A} such that  $\Vert \mathcal{A}\Vert_{\mu}<\infty$ in corresponding light tail, moderate tail or heavy tail situations , then the solution of SPDEs (\ref{whitenoisemodel})  exists in $l^2(Z^d,\mu)$, unique and continuous in time. 
\end{theorem}

\begin{proof}
Consider the classical Picards iterative method.

If we begin with $u_t^{(0)}=u(0,x)=1\in l^2(Z^d,\mu)$ and use the  notation $u_t^{(n)}=u^{(n)}(t,x)\in l^2(Z^d,\mu)$ for simplicity,
\begin{equation}
u_t^{(n+1)}=1+\int_0^t\mathcal{A}u^{(n)}(s,x)ds+\int_0^tu_t^{(n)}(s,x)dW(s,x).
\end{equation}
 
We can get a sequence of random functions with values in the space $l^2(Z^d,\mu)$. Then
\begin{equation}
<|u_t^{(k+1)}-u_t^{(k)}|^2>\leq (2T \Vert \mathcal{A} \Vert_{\mu}^2+2\sum\limits_{z\in Z^d} b^2(x-z))\int_0^t <|u_t^{(k)}-u_t^{(k-1)}|^2>ds
\end{equation}

for $k\geq 1$, $t\leq T$ and
\begin{align*}
<| u_t^{(1)}-u_t^{(0)}|^2>&=<(\int_0^tdW(s,x))^2>=\sum\limits_{z\in Z^d}b^2(x-z)t=A_1t
\end{align*} 
where $A_1=\sum\limits_{z\in Z^d}b^2(x-z)=B(0)<\infty$.

$$<\Vert u_t^{(1)}-u_t^{(0)} \Vert_{\mu}>=\sum\limits_{x\in Z^d}\mu(x)<| u_t^{(1)}-u_t^{(0)}|^2>=A_1t\sum\limits_{x\in Z^d}\mu(x)<\infty.$$ 

By induction, assume $<| u_t^{(k)}-u_t^{(k-1)}|^2>\leq \displaystyle\frac{A_1A_2^{k-1}t^{k}}{k!}$, $k\geq 2$, then we obtain 
\begin{align*}
<|u_t^{(k+1)}-u_t^{(k)}|^2>&\leq(2T \Vert \mathcal{A} \Vert_{\mu}^2+2\sum\limits_{z\in Z^d} b^2(x-z))\int_0^t<|u_t^{(k)}-u_t^{(k-1)}|^2>ds\\
&\leq A_2\int_0^t \frac{A_1A_2^{k-1}t^k}{k!} dt
\leq \frac{A_1A_2^{k}t^{k+1}}{(k+1)!}
\end{align*}
where $A_2=2T \Vert \mathcal{A} \Vert_{\mu}^2+2\sum\limits_{z\in Z^d} b^2(x-z)<\infty$. Then
\begin{align*}
<\Vert u_t^{(k+1)}-u_t^{(k)}\Vert_{\mu}>&\leq\sum\limits_{x\in Z^d}\mu(x) <|u_t^{(k+1)}-u_t^{(k)}|^2>\leq\sum\limits_{x\in Z^d}\mu(x)\frac{A_1A_2^{k}t^{k+1}}{(k+1)!}.
\end{align*}
Then
\begin{align*}
<\Vert u_t^{(m)}-u_t^{(n)}\Vert_{\mu}>&=<\Vert \sum\limits_{k=n}^{m-1}u_t^{(k+1)}-u_t^{(k)}\Vert_{\mu}>
\leq \sum\limits_{k=n}^{m-1} <\Vert u_t^{(k+1)}-u_t^{(k)}\Vert_{\mu}>\\
&\leq \sum\limits_{k=n}^{m-1} \sum\limits_{x\in Z^d}\mu(x)\frac{A_1A_2^{k}t^{k+1}}{(k+1)!}
= \sum\limits_{x\in Z^d}\mu(x)\sum\limits_{k=n}^{m-1}\frac{A_1A_2^{k}t^{k+1}}{(k+1)!}\rightarrow 0 
\end{align*}
as $m,n\rightarrow \infty$. Therefore $u_t^{(n)}$ is a Cauchy sequence in $l^2(Z^d,\mu)$. Define $u_t=\lim\limits_{n\rightarrow\infty}u_t^{(n)}$, then $u_t$ is $\mathcal{F}_t$ measurable for all $t\geq 0$, since this holds for each $u_t^{(n)}$, now we prove that $u_t$ satisfies equation (\ref{usolution}). 
For all $n$ and all $t \in[0,T]$, we have 
\begin{align*}
u_t^{(n+1)}=1+\int_0^t\mathcal{A}u^{(n)}(s,x)ds+\int_0^tu^{(n)}(s,x)dW(s,x).
\end{align*}
Now, let $n\rightarrow \infty$, then by H\"{o}lder inequality and the boundness of the $\Vert A \Vert_{\mu}$, we get that in $l^2(Z^d,\mu)$,
\begin{align*}
\int_0^t\mathcal{A}u^{(n)}(s,x)ds\rightarrow \int_0^t \mathcal{A}u(s,x)ds
\end{align*}
and 
\begin{align*}
\int_0^tu^{(n)}(s,x)dW(s,x) \rightarrow \int_0^tu(s,x)dW(s,x).
\end{align*}
We conclude that for all $t\in [0,T]$, we have 
$u(t,x)=1+\displaystyle\int_0^t\mathcal{A}u(s,x)ds+\int_0^tu(s,x)dW(s,x)$. And it exists in $l^2(Z^d,\mu)$. 
The same arguments give the uniqueness of the solution. Namely if there are two solutions $u_1(t,x)$, $u_2(t,x)$ for the equation (\ref{whitenoisemodel})
then like above 
\begin{equation*}
<|u_1(t,x)- u_2(t,x)|^2 > \leq  |2t \Vert \mathcal{A} \Vert_{\mu}^2+2B(0)| \int_0^t \langle |u_1-u_2|^2 \rangle ds
\end{equation*}
and Gronwall inequality gives $u_1 =u_2\ P-as.$

To prove the continuity of the solution $u(t,\cdot)$ as the random process in $l^2(Z^d, \mu)$, we apply similar calculations to the forth moment and prove that 
\begin{equation*}
< |u(t_1,x) - u(t_2,x) |^4 > \leq C |t_1 -t_2|^2
\end{equation*}
like in the case of the Wiener process $W(t,x)$ continuity follows now from the well-known Kolmogorov's criterion.

\end{proof}

\section{Kac-Feynman representation of the solution $u(t,x)$}
In this section, we'll prove the Kac-Feynmann representation of the It\^{o} solution $u(t,x)$ in the form of the integral over the trajectories of the random walk $x(t)$ associated to generator $\varkappa \mathcal{A}$. This representation will not be used at present but will be critical in the second part of the paper , concerning $P$-a.s. asymptotic behavior of $u(t,x), \ t \rightarrow \infty$ and calculations of $P$-a.s. Lyapunov exponents $\tilde{\gamma} = \lim_{t \rightarrow \infty} \frac{\ln u(t,x)}{t}$ in the next section. 
\begin{theorem}
Assume $\sum\limits_{z\in Z^d}b^2(z)<\infty$, $\sum\limits_{z\in Z^d}\mu(x)<\infty$, and conditions in Theorem \ref{thm for boundness of A} which guarantee $\Vert \mathcal{A}\Vert_{\mu}<\infty$ are satisfied, then the It\^{o} solution of the differential equation (\ref{whitenoisemodel}) exists and it has the representation 
\begin{align}\label{u(t,x)kac}
 u(t,x)=E_x[e^{\int_0^tdW(s,x(t-s))-\frac{tB(0)}{2}}]
\end{align}
where $x(s), s\geq 0$ is the random walk associated with the generator $\kappa\mathcal{A}$ and $E_x[\cdot]$ is the expectation over the law of $x(s),s\geq 0$ conditioned on the initial location $x$ for fixed random environment $W(\cdot,\cdot)$. 
\end{theorem}

\begin{proof}

Assume $u(t,x)=e^{\varkappa \mathcal{A}t}C(t,x)$, substitute this into equation (\ref{whitenoisemodel}),we get
\begin{align*}
C'(t,x)&=e^{-\varkappa \mathcal{A}t}\dot{W}(t,x)u(t,x)\\
C(t,x)&=1+\int_0^te^{-\varkappa \mathcal{A}s}u(s,y)dW(s,y).
\end{align*}
Thus
\begin{align*}
u(t,x) =1+\int_{0}^{t} e^{\varkappa \mathcal{A} (t-s) }u(s,y)d W(s,y)=1+\int_{0}^t \sum\limits_{y\in Z^d} p(t-s,x,y)u(s,y)dW(s, y).
\end{align*}
 
where $p(t,x,y) = P ( x(t)=y |x(0)=x)$ is the transition probability from initial state $x$ to state $y$ at time $t$. Thus, in the probability space of trajectories of the random walk $x(t)$, 
\begin{equation}\label{expression of u in expect}
u(t,x) =1+ E_x \left[ \int_{0}^t u(s, x(t-s))dW(s, x(t-s)) \right].
\end{equation}
Let us consider Kac-Feynmann type representation below:
$$
u(t,x) = E_x \left[ e^{\int_0^t d W (s, x(t-s)) -\frac{B(0)t}{2}}\right].
$$
For the fixed trajectory $x(t):\  x(0)=x,\  x(t_1)=x_1, x(t_2)=x_2, \cdots$, the integral $\displaystyle\int_0^t d W(s, x(t-s)) $ is the Winner process with mean 0 and variance $B(0)t$ . Define $V(t,x)=\int_0^t d W (s, x(t-s)) -\frac{B(0)t}{2}$, by It\'{o} formula, 
$de^{V(t,x)}=e^{V(t,x)}dW(t,x)$, thus $e^{V(t,x)}$ is a martingale and 
\begin{align*}
e^{V(t,x)}=1+E_{x}[\int_0^te^{V(s,x(t-s))}dW(s,x(t-s))].
\end{align*}

Compare this result with the equation (\ref{expression of u in expect}), we get $u(t,x) = E_x \left[ e^{\int_0^t d W (s, x(t-s)) -\frac{B(0)t}{2}}\right]$.

\end{proof}

One can prove that the solution $u(t,x)$ is the exponential martingale with respect to the law of the Wiener process $W(\cdot,\cdot)$. In particular, $<u(t,x)>=m_1(t,x)=1$. Now we will study higher moments of the solution $u(t,x)$.

\section{Equations for the Moments }
This section is devoted to the proof of the existence of the moments of all orders of the solution $u(t,x)$ to the Parabolic Anderson model (\ref{whitenoisemodel}) with correlated white noise $\dot{W}(t,x)$. First, let us derive the equations for the moments. 
For positive integer $p \geq 1$, $t \geq 0$, let  $x_1$, $x_2$, $x_3, \cdots$ be fixed points in $Z^d$ and  
denote $p^{th} $ moment by $m_p(t, x_1,x_2, \cdots, x_p)=< u(t,x_1)u(t,x_2) \cdots u(t,x_p)>.$
By Feynman-Kac formula (\ref{u(t,x)kac}) and the spatial homogeneity of the field $u(t,x)$,
\begin{align*}
m_p(t, x_1,x_2, \cdots, x_p)\leq \frac{<u^p(t,x_1)+\cdots+u^p(t,x_p)>}{p}=<u^p(t,x)>.
\end{align*}
By Lyapunov inequality, for $p\geq 1$, 
\begin{align*}
<u(t,x)>=E_x \left[ e^{\int_0^t d W (s, x(t-s)) -\frac{B(0)t}{2}}\right]\leq \left(E_x \left[ e^{p\int_0^t d W (s, x(t-s)) -\frac{B(0)t}{2}}\right]\right)^{\frac{1}{p}}.
\end{align*}
Thus
\begin{align*}
<u^p(t,x)>=E_x\left[<e^{p\int_0^t d W (s, x(t-s)) -\frac{pB(0)t}{2}}>\right]=e^{\frac{p(p-1)B(0)t}{2}}<\infty.
\end{align*}
This result implies that $m_p(t, x_1,x_2, \cdots, x_p)< \infty$.

Now we can derive the equations for $m_p(t, x_1,x_2, \cdots, x_p)$ using the It\^{o} formula. Let $f_p(t, w)=u(t, x_1)\cdots u(t, x_p) $, then 
\begin{align*}
	d f_p(t, w)= d(u(t, x_1))u(t, x_2)\cdots u(t, x_p)  + u(t, x_1) \cdots u(t, x_{n-1}) d u(t, x_p) + \sum_{i < j}  du(t, x_i) du(t, x_j)
\end{align*}
where $du(t, x_i) = \varkappa \mathcal{A} u(t, x_i)dt + u(t, x_i)dW(t, x_i)$. 
Since $dW(t, x_i)dW(t, x_j) = B(x_i - x_j)dt $, the following stochastic differential equation is obtained:
\begin{align*}
d m_p(t, x_1, \cdots, x_p)= \sum_{i=1}^p\varkappa \mathcal{A}_{x_i} m_p(t, x_1, \cdots, x_p) dt+ m_p(t, x_1, \cdots, x_p) B_p(x_1, \cdots, x_p) dt 
\end{align*}
where $\mathcal{A}_{x_i}m_p= \sum\limits_{z \neq 0, 
                  z \in Z^d }
                a(z)\big (m_p(t,x_1,\cdots,x_i+z,\cdots,x_p)-m_p(t,x_1,\cdots,x_i,\cdots,x_p)\big )$, $i=1,2,\cdots,p$ and $B_p(x_1, \cdots, x_p) = \sum\limits_{i <j} B(x_i-x_j)$. 
Finally, 
\begin{equation} \label{eq: p moment}
	\frac{d m_p(t, x_1, \cdots, x_p) }{dt} = \varkappa \left(\sum_{i=1}^p \mathcal{A}_{x_i} m_p(t, x_1, \cdots, x_p) \right) + \left( \sum_{i <j} B(x_i-x_j) \right)  m_p(t, x_1, \cdots, x_p)  
\end{equation}
with the initial condittion $m_p(0, x_1, \cdots, x_p) =1$.

In more details, the equation for first moment $m_1(t,x)= \langle u(t,x) \rangle$ is 
\begin{equation}\label{m1}
\begin{split}
  &\frac{\partial{m_1(t,x)}}{\partial{t}} =\varkappa \mathcal{A}  m_1(t,x) , \ (t,x) \in [0, \infty)\times Z^d\\
  &m_1(0,x)=1 .\\
 \end{split}
\end{equation}
Then $m_1(t,x)\equiv 1$.

For 2nd moment $m_2(t,x_1, x_2)= \langle u(t,x_1) u(t,x_2) \rangle$ , we have
\begin{equation}\label{eq: 2nd moment}
\begin{split}
  &\frac{\partial m_2(t,x_1, x_2)}{\partial t}=\varkappa (\mathcal{A}_{x_1}+\mathcal{A}_{x_2})m_2(t,x_1,x_2)+B(x_1-x_2)m_2(t,x_1, x_2)\\
  & m_2(0,x_1, x_2)= \langle u(0,x_1)u(0,x_2)\rangle=1.
 \end{split}
\end{equation}
And due to the fact that for fixed $t$, the field $u(t,x)$ is homogeneous in space, $m_2(t,x_1,x_2)=m_2(t,x_1-x_2)=m_2(t,x_1-x_2)=m_2(t,v)$, thus preceding equation is equivalent to
\begin{align}
&\frac{\partial m_2(t,v)}{\partial t}=\mathcal{H}_2m_2(t,v), \mathcal{H}_2 = 2\varkappa\mathcal{A}_v +B(v)\label{m2}\\
&m_2(0,v)=1\label{m2_initial}
\end{align}
Define $m_2(t,v)=m_{21}(t,v)+1$, then 
\begin{align}
		&\frac{\partial m_{21}(t,v)}{\partial t}=2\varkappa\mathcal{A}_v m_{21}(t,v)+B(v)m_{21}(t,v)+B(v)\label{m21}\\
		&m_{21}(0,v)=0.\label{m21_initial}
\end{align}

In summary, we get the following theorem comparable to the result in \cite{carmona1994parabolic}:
\begin{theorem} \label{thm: moments}
For each integer $p\geq 1$, each $t \geq 0 $ and $\vec{x} = (x_1,\cdots,x_p ) \in Z^{pd}$ let us set 
$$
m_p(t,\vec{x})= m_p (t, x_1, \cdots, x_p) = \left<u(t,x_1) \cdots u(t,x_p)\right>.
$$
Then these moments(Correlation functions) satisfy the following "p-particle" parabolic equation 
\begin{equation} \label{moment-corr}
  \begin{split}
    &\frac{\partial{m_p}}{\partial t} = \varkappa(\mathcal{A}_{x_1}+\cdots +\mathcal{A}_{x_p} ) m_p + \left( \sum _{i<j} B(x_i-x_j)\right)m_p\\
    & m_p(0,x) \equiv 1
  \end{split}
\end{equation}
where $\mathcal{A}_{x_i}m_p= \sum\limits_{z \neq 0, 
                  z \in Z^d }
                a(z)\big (m_p(t,x_1,\cdots,x_i+z,\cdots,x_p)-m_p(t,x_1,\cdots,x_i,\cdots,x_p)\big )$, $i=1,2,\cdots,p$. Equation (\ref{moment-corr}) can be also written as 
\begin{equation} \label{moment-corr 2}
   \begin{split}
     &\frac{\partial {m_p}}{\partial t} = \mathcal{H}_p m_p, \  \mathcal{H}_p =\varkappa(\mathcal{A}_{x_1}+\cdots +\mathcal{A}_{x_p} ) +V_p(x) \\
     & V_p(x) = \sum _{i<j}B(x_i-x_j),\ p>1
   \end{split}
\end{equation}
and Hamiltonian $ H_p$ is a classical "$p$-particle " Schr\"{o}dinger operator on the lattice $Z^{pd}$ with the binary interaction $B(x-y)$.
\end{theorem}
Equation (\ref{moment-corr 2}) will be studied in $l^2(Z^d)$. The spectral analysis of the "$p$-particle " Schr\"{o}dinger operator in $l^2(Z^d)$  or $L^2(R^d)$ plays a critical role in modern mathematical physics. The solutions of the moments equations (\ref{moment-corr 2}) have an exponential behavior as $t\rightarrow \infty$. The main purpose of this work is to investigate this asymptotic behavior. The first step is by the following result:

\begin{theorem}
For each integer $p\geq 1$, the limit 
\begin{align*}
\lim\limits_{t\rightarrow\infty}\frac{1}{t}\ln m_p(t,x)
\end{align*}
is independent of $\boldsymbol{x}=(x_1,x_2,\cdots,x_d)$. This limit is called the $p$-th (moment) Lyapunov exponent of the solution $u(t,x)$ and it is denoted by $\gamma_p(\varkappa)$. It is equal to the supremum of the spectrum of the multiparticle Schr\"{o}dinger operator $\mathcal{H}_p$ appearing in the moment equation. 
\end{theorem}

The proof of this result is standard and the similar proof for the local Laplacian operator $\Delta$ is given in full details at \cite{carmona1994parabolic}. We will concentrate on the second moment $m_2(t,v)$, $v=x_1-x_2$ corresponding to the Hamiltonian $\mathcal{H}_2=2\varkappa \mathcal{A}+B$ and Lyapunov exponent $\gamma_2(\varkappa)=\max\{\lambda\in Sp(\mathcal{H}_2)\}$ where $Sp(\mathcal{H}_2)$ stands for the spectrum of $\mathcal{H}_2$.And the existence of $p$-th moment Lyapunov exponent depends significantly on the spectral analysis of the semigroup $T_{+}=\exp\{2\varkappa t\mathcal{A}\}$ and its kernel $p(t,x,y)$. 

\section{Transition Probability $p(t,x,y)$}
Let us find out the transition probability $p(t,x,y)=P(x(t)=y|x(0)=x)$ of the random walk $x(t)$ associated with $2\varkappa\mathcal{A}$. It is well known that $p(t,x,y)$ is the fundamental solution of the following parabolic problem
\begin{align}\label{p(t,x,y)equation}
\frac{dp}{dt}&=2\varkappa\mathcal{A} p, t>0\\
p(0,x,y)&=\delta_y(x)\nonumber.
\end{align}

Applying the Fourier transform technique at both sides of equation(\ref{p(t,x,y)equation}), where Fourier transform of $p$ is defined as $\hat{p}(t,x,k)=\sum\limits_{y\in Z^d}p(t,x,y)e^{i k  y}$, we will get the solution to (\ref{p(t,x,y)equation}):

\begin{equation}
p(t,x,y)=\frac{1}{(2\pi)^d}\int_{T^d}e^{2\varkappa \hat{\mathcal{A}}(k)t}e^{i k (x-y)} dk, \,\,\, k\in T^d=[-\pi,\pi]^d
\end{equation}

where $\hat{\mathcal{A}}(k)=\sum\limits_{z\neq 0}(cos(k,z)-1)a(z)=\hat{a}(k)-1$ is the Fourier symbol of the operator $\mathcal{A}$.

Thus
\begin{equation} \label{p(t,0,0)}
p(t,x,x)=p(t,0,0)=\displaystyle\frac{1}{(2\pi)^d}\int_{T^d}e^{2\varkappa \hat{\mathcal{A}}(k)t}dk,\,\,\,k\in T^d.
\end{equation}

In the future, we'll study the phase transition with represent to $\varkappa$ of the top eigenvalue for the Hamiltonian $\mathcal{H}_2 = 2 \varkappa \mathcal{A}+B(x), B(x) \geq 0, \  x \in Z^d$. In this situation it is convenient to introduce the standard diffusivity $\varkappa_0 = \frac{1}{2}$, i.e., $2 \varkappa_0 =1$. Let's denote in this case
$$p_0(t,x,y) = \frac{1}{(2 \pi)^{d}} \int_{T^d} e^{t \hat{\mathcal{A}}(k)+i k(x-y) } dk.$$
Then (for arbitrary $\varkappa$)
$$
p(t,x,y) = p_0(2 \varkappa t ,x,y).
$$
Due to Central Limit Theorem in the case of finite second moment $\sigma^2 = \sum\limits_{z \in Z^d} |z|^2a(z)$,
$$
p_0(t,x,y) \sim \frac{exp \left[- \frac{\Pi ^{-1} \{x-y,x-y \}}{2} \right]}{(2 \pi t)^{d/2} \sqrt{det \Pi}} 
$$
if $t \rightarrow  \infty $, $ |x-y| = \underline{\underline{O}}(\sqrt{t})$. Here $\Pi$ is the correlation matrix of the random walk $x(t)$ in the case $2 \varkappa_0 = 1$:
$$
\Pi = - \left[ - \frac{\partial^2 \hat{A}(0)}{\partial \theta_i \partial \theta_j}  \right].
$$ 
Because of the general assumption that $\Pi$ is positive definite, $det(\Pi)>0$. It means that for appropriate constant $C^+$ depending only on $\hat{a}(k)$ and all $t \geq 1$ we have estimate
\begin{equation} \label{eq:p_0}
p_0(t,x,x) \leq \frac{C^+}{t^{d/2}}, \ t \geq 1 .
\end{equation}

The following well known result can be proved based on the fact of (\ref{eq:p_0}).

\begin{theorem}
If $a(z)$ satisfies light tail condition (\ref{light tails condition}) or moderate tail condition (\ref{moderate tails condition}), then $\sigma^2 = \sum\limits_{z \neq 0 } |z|^2a(z) < \infty$, and the random walk $x(t)$ associated with $2\kappa\mathcal{A}$ has transition probability $p(t,x,x) \sim \displaystyle\frac{C}{t^{d/2}}$, then $x(t)$ is transient when $d\geq 3$ and $x(t)$ is recurrent when $d=1,2$. 
\end{theorem}

\begin{theorem}
Assume that $\sigma^2 = \sum\limits_{z \neq 0 } |z|^2a(z) < \infty$ but $a(z)$ satisfies heavy tail condition (\ref{heavy tails condition}), then
$$
p_0(t,x,x) =p_0(t,0,0) \sim \frac{C^+}{t^{d/ \alpha}}, \ 0< \alpha <2.
$$
\end{theorem}

The following theorem gives the summary of the results:
\begin{theorem}
Denote the random walk generated by the operator $2\varkappa \mathcal{A}$ on $Z^d$ as $x(t)$.
\begin{enumerate}
\item If $d\geq 3$, then the random walk $x(t)$ is transient without any additional technical conditions. 
\item If $d=2$, then under regularity condition (\ref{heavy tails condition}), the random walk $x(t)$ is transient for $0<\alpha<2$. If  $\sigma^2=\sum\limits_{z\neq0}|z|^2a(z)<\infty$, $\alpha=2$, the random walk $x(t)$ is recurrent. 
\item If $d=1$ and $\alpha\leq 1$, the random walk $x(t)$ is transient, otherwise $x(t)$ is recurrent.

\end{enumerate}
\end{theorem}

\section{ Spectral analysis of the operator $\mathcal{H}_2 = 2  \varkappa \mathcal{A}+B(x)$ and the problem of intermittency}

The first moment of the field $u(t,x)$ equals to $1$ identically. The fluctuations of this field for $t \rightarrow \infty$ depends on higher moments and first of all, on the second moment $m_2(t,x,y)=m_2(t,0, y-x)=m_2(t,v)$. Recall $m_2(t,v)$ satisfies equation (\ref{m2}) and initial condition (\ref{m2_initial}): 
\begin{align}\label{spectralm2equ}
\frac{\partial m_2(t,v)}{\partial t}&=\mathcal{H}_2m_2(t,v)\\
m_2(0,v)&=1
\end{align}
where $\mathcal{H}_2=2\varkappa\mathcal{A}_v+B(v)$.  Operator $\mathcal{H}_2$ is bounded and $B(x)\rightarrow 0$ as $x\rightarrow \infty$. The spectrum of $\mathcal{H}_2$ consists of two parts. The first part is the essential spectrum of $\mathcal{H}_2$ denoted by $Sp_{ess}\mathcal{H}_2=[2\varkappa c,0]$ where $c=\min\limits_{k\in T^d}(\hat{a}(k)-1)$, and $Sp_{ess}\mathcal{H}_2$  equals the spectrum of $\mathcal{H}_2=2\varkappa\mathcal{A}$. The second part is at most countable discrete spectrum denoted by $Sp_d(\mathcal{H}_2)$ with possible accumulation points $2\varkappa c$ and $0$. Denote $\lambda_0(\mathcal{H}_2)$ as the maximum point of $Sp(\mathcal{H}_2)$. If $\lambda_0(\mathcal{H}_2)>0$, then the positive discrete spectrum is not empty and $\lambda_0(\mathcal{H}_2)$ is the top eigenvalue of $\mathcal{H}_2$. It is necessarily simple and corresponding eigenfunction $\psi_0(x)$ is strictly positive. 

If $\lambda_0(\mathcal{H}_2)=0$, then the second moment is bounded and the fluctuations of $u(t,x)$ are bounded, this is the regular case. If $\lambda_0(\mathcal{H}_2)>0$, $m_2(t,v)$, $v=x_1-x_2$, tends to $\infty$ exponentially, then the intermittency phenomenon arises. Zeldovich et al. \cite{zel1987intermittency}, G\"{a}rtner and Molchanov \cite{gartner1990parabolic} developed the mathematical theory of intermittency with applications to hydrodynamics, magnetic and temperature fields in the turbulent flow etc.  For example, the solar magnetic field has an intermittent structure because more than 99\% of the magnetic energy concentrates on less than 1\% of the surface area.  In this situation, the main contributions to the magnetic field are from the very high and very sparse peaks (black spots). 

We will discuss in this section that fundamental phase transition from the regular case to the intermittent structure of the field $u(t,x)$ and other spectral bifurcation depending on the rate of the jumps $\{a(z),z\in Z^d\}$ and the correlation function $B(x), x\in Z^d$ which is a positive definite potential in equation for the second moment (\ref{spectralm2equ}). Part of the results in this paper can be applicable to the case with the general potential $V(x), V(x)\rightarrow 0, x\rightarrow \infty$.

\subsection{Spectral bifurcation caused by $\hat{a}(k)$}

Now we will present the spectral bifurcation depending on the jump distribution $\{a(z),z\neq 0\}$. Recall the non-local Laplacian operator $2\varkappa\mathcal{A}$ has Fourier transform representation $2\varkappa\hat{\mathcal{A}}(k)=2\varkappa(\hat{a}(k)-1)$ in $L^2(T^d,dk)$.
\begin{theorem}
If $\hat{a}(k)$ is analytic, that is $a(z)\leq ce^{-\eta|z|}$ (light tail), $c,\eta>0$, then the spectral measure of $2\varkappa\mathcal{A}$ is absolutely continuous (a.c.).
\end{theorem}

\begin{proof}
Let $\Gamma_{\lambda}=\{k\in Z^d:2\varkappa\hat{\mathcal{A}}(k)<\lambda\}$. then $E_{\lambda}f(k)=I_{\Gamma_{\lambda}}(k)f(k)$ is the family of the spectral projections associated to self-adjoint operator $2\varkappa\mathcal{A}$. To prove the absolute continuity of the spectral measure, it is sufficient to prove that 
\begin{align*}
m(k\in \Gamma_{\lambda+d\lambda}-\Gamma_{\lambda})=\rho(\lambda)d\lambda, \,\,\rho(\lambda)\in L^1(Sp(2\varkappa\mathcal{A}))
\end{align*}

where $m(k\in \Gamma_{\lambda+d\lambda}-\Gamma_{\lambda})$ is the distribution of $2\varkappa\hat{\mathcal{A}}(k)$.
Due to analyticity of $\hat{\mathcal{A}}(k)$, however,
\begin{align*}
m(k\in \Gamma_{\lambda+d \lambda}-\Gamma_{\lambda})=d \lambda\int_{k:2\varkappa \hat{\mathcal{A}}(k)=\lambda}\frac{dk}{2\varkappa|\nabla\hat{\mathcal{A}}(k)|}=d \lambda\rho(\lambda)
\end{align*}

and the integral for $\rho(\lambda)$ is well defined for all $\lambda$ except finitely many points (values of $2\varkappa\hat{\mathcal{A}}(k)$ in the critical points where $\nabla \hat{\mathcal{A}}(k)=0$). 
\end{proof}
With the similar proof, we can get the following theorem:
\begin{theorem}
If $\hat{\mathcal{A}}(k)\in C^2(T^d)$ and $\nabla\hat{\mathcal{A}}(k)=0$ only in the finite set of points $k_1,k_2,\cdots,k_n$, $\Vert Hess(\hat{\mathcal{A}}(k_j))\Vert<\infty$, $j=1,2,\cdots n$ and $det(Hess(\hat{\mathcal{A}}(k_j))\neq 0$, that is $\hat{\mathcal{A}}(k)$ is a Morse type function, then the spectral measure of $2\varkappa \mathcal{A}$ is absolutely continuous. 
\end{theorem}
If $\hat{\mathcal{A}}(k)$ is not analytic but still $C^{\infty}$ function, the situation is different. 
\begin{theorem}
If $\hat{\mathcal{A}}(k)$ is not analytic but belongs only to $C^{\infty}$-class, then one can find jumps law $\{a(z),z\in Z^d\}$ with assumptions of symmetry and aperiodicity such that $2\varkappa\hat{\mathcal{A}}(k)=c_0$, $c_0\in S_0\subset T^d$ where $S_0$ is an open ball on the torus $T^d$. In this case, $\mathcal{H}_2$ has positive eigenvalue $\lambda_0\in Sp(2\varkappa\mathcal{A})$ of the infinite multiplicity embedded into $Sp_{ess}\mathcal{H}_2$. That is, the spectral measure is not pure absolutely continuous. 
\end{theorem}
\begin{proof}
We can prove this theorem the by following particular example.
Suppose $d =1$, $k\in [- \pi, \pi] = T^1$, there exists $\hat{a}(k)=\sum\limits_{x\in Z^1}\cos(kx) a(x)$ such that  $\hat{a}(k)$ is constant when $ k \in [\alpha,\beta]$(see \cite{molchanov2007quasicumulants}).Then operator $2\varkappa\mathcal{A}$ has infinite dimensional eigenspace embedded into $Sp_{ess}\mathcal{H}_2$. Namely, if $\hat{f}(k)$ is supported on $[\alpha, \beta]$, then 
$
\varkappa \left( \widehat{\mathcal{A}f}\right)(k)  = \varkappa(\hat{a}(k)-1) \hat{f}(k) =\varkappa h \hat{f}(k)
$
where $h=\hat{a}(k)-1$. 
\end{proof}
The function $\hat{a}(k)$ may be not smooth but we can improve situation if consider the convolution $\hat{a}(k) \ast b_0(k)$, where $b_0(k)\in  C^{\infty}$ is positive definite and compactly supported on $[- \varepsilon, \varepsilon]$ and $\varepsilon < \displaystyle\frac{|\beta-\alpha|}{3}$. To construct such $b_0(k)$, it is sufficient to take any real symmetric function $\phi(k)$ supported on $[-\frac{\epsilon}{2},\frac{\epsilon}{2}]$ and put $b_0(k)=\phi(k)\ast\phi(-k)=\phi(k)\ast\phi(k)$, then $B_0(x)=\displaystyle\frac{1}{(2\pi)^d}\int_{T^d}\phi(k)\ast\phi(-k)e^{-ikx}dk=\Phi^2(x)\geq 0$ where $\Phi(x)=\displaystyle\frac{1}{(2\pi)^d}\int \phi(k)e^{-ikx}dk$. That is, we guarantee the positive definite property of the new probability jumping kernel with $\{a(x)B_0(x),x\in Z^d\}$ after the normalization. At the same time, $\hat{a}(k) \ast b_0(k)$ is constant on some subinterval $[\alpha,\beta]$.

\section{Spectral bifurcation in recurrent case and transient case}
This section is devoted to the positive discrete spectrum, that is, the problem of existence of $\lambda_0(\mathcal{H}_2)>0$. Let $N_0(B)=\sharp\{\lambda_j:\lambda_j(\mathcal{H}_2)>0\}$ be the number of positive eigenvalues of the operator $\mathcal{H}_2$ with the potential $B$. Our goal is to estimate $N_0(B)$. We will discuss $N_0(B)$ for transient case and recurrent case. Let us note that $\Vert \mathcal{H}_2\Vert\leq 2\varkappa\max\limits_{k\in T^d}|\hat{a}(k)-1|\leq 4\varkappa$. 

\subsection{Recurrent case}
The following theorem is the generalization of the similar result in (\cite{carmona1994parabolic}).
\begin{theorem}\label{recurrentthm1}
Assume that $B(x)=\delta_0(x)$, that is $W(t,x)=w(t,x), x \in Z^d$ are i.i.d Wiener processes, then $\lambda_0(\mathcal{H}_2)>0$ for any $\varkappa >0$ if and only if $\displaystyle\int_{T^d} \frac{dk}{1-\hat{a}(k)}=\infty$. i.e., random walk $x(t)$ associated with $2\varkappa\mathcal{A}$ is recurrent. 
\end{theorem}
\begin{proof}
To solve the equation $2\varkappa A\psi+\delta(x)\psi=\lambda\psi$ where $\psi$ is the eigenfunction corresponding to the eigenvalue $\lambda$, we can apply Fourier transform:
\begin{align}
2\varkappa \hat{\mathcal{A}}(k)\hat{\psi}+\psi(0)&=\lambda\hat{\psi}\\
\hat{\psi}(k)&=\frac{\psi(0)}{\lambda+2\varkappa(1-\hat{a}(k))}\\
1&=\frac{1}{(2\pi)^d}\int_{T^d}\frac{dk}{\lambda+2\varkappa(1-\hat{a}(k))}=I(\lambda). \label{I(lambda)equ}
\end{align}
The $I(\lambda)$ in the right hand of the above equation is a monotone decreasing function of $\lambda$ and $I(\lambda)\rightarrow \displaystyle\frac{1}{(2\pi)^d}\int_{T^d}\frac{dk}{2\varkappa(1-\hat{a}(k))}=I(0)$. 

It is well known that the classification of the random walk depends on the Green function $G_{\lambda}(t,x,y)$ for $p(t,x,y)$ defined as 
\begin{align*}
G_{\lambda}(t,x,y)=\int_0^te^{-\lambda u}p(u,x,y)du, \,\,\,\, \lambda \geq 0.
\end{align*}
Recall that a random walk is called recurrent if $G_0(0,0)=\infty$.
For the process $x(t)$ with generator $2\varkappa \mathcal{A}$ , 
\begin{align*}
G_{0}(0,0)=\displaystyle\frac{1}{(2 \pi)^d} \int_{T^d} \frac{dk}{2\varkappa(1 -\hat{a}(k))} =I(0).
\end{align*}
When the process $x(t)$ is recurrent, $G_0(0,0)=I(0)=\infty$, then there exists unique solution $\lambda_0=\lambda(\varkappa)>0$. 
\end{proof}

If $\varkappa=\varkappa_{cr}$, the answer depends on the dimension. 
\begin{theorem}\label{B(x_0)>0positive_general}
If the process $x(t)$ associated with $2\varkappa\mathcal{A}$ is recurrent, then $\lambda_0(\mathcal{H}_2)>0$ for nonnegative potential $B(x)\geq 0$ and $B(0)=\max\limits_{x\in Z^d}B(x) >0$.
\end{theorem}
\begin{proof}
By variation principle, to find $\lambda_0(\mathcal{H}_2) = \max\limits_{\substack{\psi \in l^2(Z^d)}:\| \psi\|=1} (\mathcal{H}_2 \psi, \psi) = (\mathcal{H}_2 \psi_0,\psi_0 )$, where $\psi_0$ is the eigenfunction corresponding to the top eigenvalue,  since  $\lambda_0(\mathcal{H}_2)\geq \lambda_0(\tilde{H})$, it is sufficient to estimate $\lambda_0(\mathcal{H}_2)$ by the eigenvalue of $\tilde{H} =  2\varkappa \mathcal{A} + B(x) \delta_{0}(x)$ which satisfies the following equation  
\begin{equation}\label{eq:I(lambda)_recurrent}
\frac{1}{B(0)} = \frac{1}{(2 \pi)^d} \int_{T^d} \frac{dk}{\lambda + 2 \varkappa (1-\hat{a}(k))} = I(\lambda).
\end{equation}
In the recurrent case, due to $G_0(0,0)=\infty$, $I(\lambda) \rightarrow \infty$ as $\lambda \downarrow 0 $ and equation (\ref{eq:I(lambda)_recurrent}) has unique solution, thus there always exists  $\lambda_0(\tilde{H}) >0$. 

\end{proof}
The following corollary follows Theorem \ref{B(x_0)>0positive_general}:

\begin{corollary}
Assume that $B(x)\geq 0, x \in Z^d$, $B(0)>0$, and $x(t)$ is recurrent, then the second moment of $u(t,x)$ will grow exponentially.
\end{corollary}

\subsection{Transient Case}

The paper by S. Molchanov and B. Vainberg \cite{molchanov2010general} contains the following the Zwickel-Lieb-Rozenblum (ZLR) type estimate for the number of operator $\mathcal{H}_0$'s positive eigenvalues denoted by $N_0(V)=\sharp\{\lambda_i:\lambda_i(\mathcal{H}_0)>0\}$ where $\mathcal{H}_0=2\varkappa\mathcal{A}+V(x)$ with arbitrary general potential $V$, for any arbitrary constant $\sigma>0$ and in transient case,
\begin{align}
N_0(V)\leq\sharp\{x_j: |V(x_j)|\geq 4\varkappa\geq \Vert 2\varkappa\mathcal{A}\Vert\}+\frac{1}{C(\sigma)}\sum\limits_{x:0<|V(x)|<4\varkappa}V(x)\int_{\frac{\sigma}{|V(x)|}}p(t,x,x)dt
\end{align}
where $p(t,x,x)$ is the transition probability for the random walk with the generator $2\varkappa\mathcal{A}$ and $C(\sigma)=e^{-\sigma}\int_0^{\infty}\frac{ze^{-z}}{z+\sigma}dz$. Then for $\sigma=2$, with transition probability $p_0(,t,x,y)$ that is associated with the generator $\mathcal{A}$,
\begin{align}\label{generalized CLR estimate}
N_0(V)\leq \sharp\{x_j: |V(x_j)|\geq 4\varkappa\geq \Vert 2\varkappa\mathcal{A}\Vert\}+\frac{1}{2\varkappa C(2)}\sum\limits_{x:|V(x)|\leq 4\varkappa}V(x)\int_{\frac{4\varkappa}{|V(x)|}}p_0(s,0,0)ds.
\end{align}
This estimate will be useful in the future.  

Our goal is to prove that in the transient case, there is an important phase transient (spectral bifurcation). One can find constant $\varkappa_0(\hat{a})$ depending only on the distribution of the jumps (or  corresponding characteristic function $\hat{a}$) such that for $\varkappa < \varkappa_0$ there is at least one positive eigenvalue , i.e., $N_0(B) \geq 1$, but $N_0(B) =0$ for $\varkappa > \varkappa_0$. It leads to phase transition in the asymptotic behavior of the solution $u(t,x)$ of the Anderson parabolic problem. The second moment $m_2$ of the solution is growing exponentially for small enough $\varkappa$, i.e., $\varkappa < \varkappa_0)$ and bounded for the large enough $\varkappa$, i.e., $\varkappa > \varkappa_0)$.  The analysis of simpler case when the generator associated with the random walk is the local lattice Laplacian and the potential $\xi_t(x)$ are independent white noises has been done in \cite{carmona1994parabolic}. However, in our case, the generator of random walk is nonlocal Laplacian and the potential is more complicated.

\begin{theorem}
Assume that $B(x)=\delta_0(x)$, that is $W(t,x)=w(t,x), x \in Z^d$ are i.i.d Wiener processes, then if $\displaystyle\int_{T^d}\frac{dk}{1-\hat{a}(k)}<\infty$, i.e., the random walk $x(t)$ is transient, then there exists a critical value $\varkappa_{cr}$ such that $\lambda_0(\mathcal{H}_2)>0$ for $\varkappa<\varkappa_{cr}$ and $\lambda_0(\mathcal{H}_2)=0$ for $\varkappa>\varkappa_{cr}$. 
\end{theorem}
\begin{proof}
The proof is analogy to Theorem \ref{recurrentthm1}.  Recall that a random walk is called transient if $G_0(0,0) < \infty$, for the process $x(t)$ with generator $2\varkappa \mathcal{A}$ , 
\begin{align*}
G_{0}(0,0)=\displaystyle\frac{1}{(2 \pi)^d} \int_{T^d} \frac{dk}{2\varkappa(1 -\hat{a}(k))} =I(0).
\end{align*}
When the process $x(t)$ is transient, due to $G_0(0,0)<\infty$, $I(\lambda) < \infty$ as $\lambda \downarrow 0 $ and equation (\ref{I(lambda)equ}) has unique solution if  $\varkappa\leq \varkappa_{cr}$ where $\varkappa_{cr}=\displaystyle\frac{1}{2 (2\pi)^d}\int_{T^d}\frac{dk}{1-\hat{a}(k)}$.  Thus, there exists a unique simple positive eigenvalue $\lambda_0>0$ with positive eigenfunction $\psi_0(x)>0$ if $\varkappa< \varkappa_{cr}$. And the lack of the solution of equation (\ref{I(lambda)equ})when $\varkappa>\varkappa_{cr}$ tells that there is no positive eigenvalue $\lambda_0>0$ if $\varkappa>\varkappa_{cr}$.

\end{proof}

\begin{theorem}\label{B(x_0)>0positive_transient}
If the process $x(t)$ is transient, the potential $B(x)\geq 0$ and $B(0)=\max\limits_{x\in Z^d}B(x) >0$, then $\lambda_0(\mathcal{H}_2) >0$ exists for small enough $\varkappa$. 
\end{theorem}
\begin{proof}
By variation principle, to find $\lambda_0(\mathcal{H}_2) = \max\limits_{\substack{\psi \in l^2(Z^d)}:\| \psi\|=1} (\mathcal{H}_2 \psi, \psi) = (\mathcal{H}_2 \psi_0,\psi_0 )$, where $\psi_0$ is the eigenfunction corresponding to the top eigenvalue,  since  $\lambda_0(\mathcal{H}_2)\geq \lambda_0(\tilde{H})$, it is sufficient to estimate $\lambda_0(\mathcal{H}_2)$ by the eigenvalue of $\tilde{H} =  2\varkappa \mathcal{A} + B(x) \delta_{0}(x)$ which satisfies the following equation  
\begin{equation}\label{eq:I(lambda)_transient}
\frac{1}{B(0)} = \frac{1}{(2 \pi)^d} \int_{T^d} \frac{dk}{\lambda + 2 \varkappa (1-\hat{a}(k))} = I(\lambda).
\end{equation}
In the transient case, due to $G_0(0,0)<\infty$, $I(\lambda) < \infty$ as $\lambda \downarrow 0 $ and equation (\ref{eq:I(lambda)_transient}) has unique solution if $I(0)> \frac{1}{B(0)}$, otherwise there is no positive eigenvalue, where $I(0)=\frac{1}{(2\pi)^d}\int_{T^d}\frac{dk}{2\kappa(1-\hat{a}(k))}$. The condition  $I(0)> \frac{1}{B(0)}$ is equivalent to $\varkappa\leq \frac{B(0)}{2 (2\pi)^d}\int_{T^d}\frac{dk}{1-\hat{a}(k)}$.  Thus, there exists a unique simple positive eigenvalue of $\tilde{H}$ with positive eigenfunction $\psi_0(x)>0$ if $\varkappa\leq \frac{B(0)}{2 (2\pi)^d}\int_{T^d}\frac{dk}{1-\hat{a}(k)}$ when the random walk is transient. 

\end{proof}

\begin{theorem} \label{th: thereom 1}
If $d \geq 3$, $\sum\limits_{z \in Z^d} |z|^2 a(z) < \infty$,  $\sum\limits_{x \in Z^d} |B(x)|^{d/2} < \infty$, $B\geq 0$, then there exists $\varkappa_{cr}$, if $\varkappa >\varkappa_{cr}$, $N_0(B)=0$, if $\varkappa < \varkappa_{cr}$, $N_0(B) \geq 1  $.
\end{theorem}

\begin{proof}
  Since $B(x) \in l^1(Z^d)$, $B(x)\rightarrow 0$ as $|x|\rightarrow\infty$. Then there exists a constant $q>0$ such that $\displaystyle\frac{2 \varkappa q}{B(x)} \geq 1$ for all $x\in Z^d$. Thus $\sharp\{x_j: |B(x_j)|\geq 4\varkappa\geq \Vert 2\varkappa\mathcal{A}\Vert_{\mu}\}<\infty$. When $d\geq 3$ and $\sum\limits_{z \in Z^d} |z|^2 a(z) < \infty$, the random walk is transient and $p(t,x,y)$ satisfies (\ref{eq:p_0}), apply ZLR estimate in (\ref{generalized CLR estimate}), we get
\begin{align}
N_0(B)  &\leq  \sharp\{x_j: |B(x_j)|\geq 4\varkappa\geq \Vert 2\varkappa\mathcal{A}\Vert_{\mu} \}+\frac{1}{2 \varkappa C(q)} \sum_{x \in Z^d} |B(x)| \int_{\frac{2 \varkappa q}{B(x)} } \frac{C_+}{s^{d/2}} ds\\
&= \sharp\{x_j: |B(x_j)|\geq 4\varkappa\geq \Vert 2\varkappa\mathcal{A}\Vert_{\mu}\}+C_1   \frac{1}{{\varkappa^{d/2}}}\sum_{x \in Z^d}{|B(x)|^{d/2}} <\infty .\label{eq: ZLR estimate}
\end{align} 
Let's stress that constant $C_1$ in the Zwickel-Lieb-Rozenblum (ZLR) type estimate (\ref{eq: ZLR estimate}) depends only on dimension $d$ and the function $\hat{\mathcal{A}}(k)$ (i.e., the distribution of jumps $a(\cdot)$ and $ B(0)$).

This estimate (\ref{eq: ZLR estimate}) proves the existence of the phase transition (spectral bifurcation): there exists a critical value $ (\varkappa_{cr})$ such that $N_0(B)=0$ for $\varkappa > \varkappa_{cr}$ and  $N_0(B)\geq 1$ for $\varkappa < \varkappa_{cr}$.

\end{proof}

\begin{remark}
The classical ZLR result concerns the $N_0(V)$ for the continuous Schr\"odinger operator $\mathcal{H} = - \Delta +\sigma V(x)$ on $L^2(R^d, dx), \ d \geq 3$. The Brownian motion associated to the Laplacian $\Delta$ is transient exactly for $d \geq 3$.  The classical ZLR estimate gives the upper bound of the number of negative eigenvalue of operator $\mathcal{H}$: $N_0(V) = \sharp  \{\lambda _i(\mathcal{H}) <0 \}$, that is
\begin{equation}
  N_0(\sigma V) \leq C_0(d) \sigma^{d/2}\int_{R^d} |V(x)|^{d/2} dx 
\end{equation}

and for large $\sigma$, it gives the right order of $N_0(\sigma V)$, which is so called quasi-classical estimate. In the lattice case, where we can put $\sigma=\frac{1}{\varkappa}$, this is not true. If, say $V(x)=\mathbb{I}_{\Gamma}(x)$, $\Gamma$ is a finite set, then for small $\varkappa$, $N_0(V)=card(\Gamma)$ that remains bounded.
\end{remark}

Let's stress that for  $d \geq 3$, any random walk on $Z^d$ is transient independently on the existence of the second moment. If $\sum\limits_{z\in Z^d}z^2a(z)=\infty$, the random walk can be also transient in any dimension and belongs to the domain of attraction of the stable and symmetric law with parameter $\alpha$ when  under some regularity conditions see \cite{agbor2015global}\cite{gartner1990parabolic}. If $d=1$, then $\alpha<1$. If $d\geq 2$, then $\alpha<2$. 

In this case, the ZLR estimate has the following form:
\begin{equation}
N_0( B) \leq \sharp\{x_j: |B(x_j)|\geq 4\varkappa\}+C_1 \sum_{x \in Z^d} \frac{|B(x)|^{d/\alpha}}{\varkappa^{d/2}}.
\end{equation}
\begin{theorem} \label{th: thereom 2}
If the random walk $x(t)$ has the limiting stable law with parameter $0< \alpha <1$ for $d=1$  and $0< \alpha <2$ for $d \geq 2 $, i.e., $\displaystyle\frac{x(t)}{t^{d/\alpha}}$ converges weakly to a stable distribution $S t_{\alpha}$ with parameter $\alpha$ in distribution, and $\sum\limits_{x \in Z^d} |B(x)|^{d/\alpha} < \infty$, then there exists $\varkappa_{cr} >0$ such that for small $\varkappa<\varkappa_{cr}$, $\ N_0(B) \geq 1$ and for $\varkappa> \varkappa_{cr}$, $N_0(B)=0$. 
\end{theorem}

Theorems \ref{th: thereom 1}, \ref{th: thereom 2} demonstrate in the transient case, the bifurcation from the regular behavior of the field $u(t,x)$ to intermittent behavior.

\begin{corollary}
In the situation of the Theorem \ref{th: thereom 1} and \ref{th: thereom 2}, for all sufficiently large $\varkappa$, there is no positive discrete spectrum, that is $N_0(B)=0$.
\end{corollary}

We have discussed the cases when  $\sum\limits_{x \in Z^d} |B(x)|^{d/2} < \infty$ and when $\sum\limits_{x \in Z^d} |B(x)|^{d/\alpha} < \infty$ with $d=1$, $\alpha<1$ or $d\geq2$, $0<\alpha<2$. The question is that what if the above series are divergent?  The following Theorem \ref{thm:sum B^beta inf} gives answer for the situation when $\sum\limits_{x\in Z^d}|B(x)|^{\beta}=\infty$ for any arbitrary large $\beta$ in a particular case.

\begin{theorem}\label{thm:sum B^beta inf}
Let $\mathcal{A}f(x)=\Delta f(x)=\sum\limits_{x':|x'-x|=1}(f(x')-f(x))$, that is, the process $x(t)$ associated with the generator $\mathcal{A}$ is the classical random walk with jumps to the nearest neighbors. $d\geq 3$ and $B(x)=\sum\limits_{n=1}\sigma_n\delta(x-x_n)$, where $\sigma_n \rightarrow 0$ very slowly (say $\sigma_n=\frac{\delta}{\ln(n+1)}$) and the points $\{x_n,n\geq 1\}$ are very sparse, say $x_n=(2^n,0,0,\cdots,0)$, $n\geq 1$. Then $N_0(B)=0$ for sufficiently small $\delta>0$.
\end{theorem}

\begin{proof}
We'll give the proof for particular case presented above, but it will be clear that the result is very general. Since $B(x) \geq 0$ and $B(x) \rightarrow 0$, the spectrum  of $\mathcal{H}_2$  for $\varkappa > 0$ is discrete. Assume that it is non-empty and $\lambda_0(\mathcal{H}_2)$ is the maximum postitive eigenvalue with positive eigenfunction $\psi_0(x)$, then Fourier transform for  $\psi_0(x_n),  \  n \geq 1$ is 
$$
-\hat{\Delta} (k) \hat{\psi}(k) + \sum_{m =1}^{\infty} \sigma_m \psi_0(x_m) e^{i(k,m)} = \lambda_0 \hat{\psi}_0(k).
$$
i.e., 
$$
\hat{\psi}_0(k) = \sum_{m =1}^{\infty} \sigma_m \psi_0(x_m) \frac{e^{i(k,m)}}{\lambda_0 +  \hat{\Delta} (k)},\,\,\,\, \hat{\Delta} (k)=\sum\limits_{j=1}^d(\cos k_j-1)
$$
and (for inverse Fourier transform)
$$
\psi_0(x_n) =  \sum_{m =1}^{\infty} \sigma_m \psi_0(x_m) G_{\lambda_0}(x_n,x_m),  \ G_{\lambda_0}(0,z)= \frac{1}{(2 \pi)^d}  \int_{0}^{\infty} \frac{e^{i(k,z)}dk}{\lambda_0 + \hat{\Delta}(k)}
.$$
Finally, 
\begin{equation} \label{psi_0}
\psi_0(x_n) = \sum_{m : m \neq n}  \frac{\sigma_m \psi_0(x_m) G_{\lambda_0}(x_n,x_m)}{1- \sigma_n G_{\lambda_0}(x_n, x_n)}
\end{equation} 

But for $d \geq 3$
$$
G_{\lambda_0}(x,y) = \int_{0}^{\infty} e^{- \lambda_0 t}p(t,x,y)dt \leq G_0(x,y) \leq \frac{C(d)}{1 +|x-y|^{d-2}}.
$$
Let's consider now the homogeneous system (\ref{psi_0}) in the Banach space $L^\infty (Z^d)$ with the norm
$$
\Vert f(x)\Vert_{\infty} = \max_{x \in Z^d} |f(x)|.
$$
Let 
\begin{equation*}
\tau_{m,n} = \begin{cases}
\frac{\sigma_m \psi_0(x_m) G_{\lambda_0}(x_n,x_m)}{1- \sigma_n G_{\lambda_0}(x_0, x_0)},\ m \neq n\\
0 , \ m = n
\end{cases}
\end{equation*}
form the matrix of (\ref{psi_0}). Assume that $|x_n - x_m| \geq |2^n - 2^m|$ for $m \neq n$, then 
$$
\Vert \tau \Vert_{\infty} \leq \max_{n} \sum_{m: m \neq n} |\tau_{m,n}| \leq C_1(d) \max_{m} \sigma_m
$$
and for sufficiently small $\max\limits_{m} \sigma_m$ we'll get $\Vert \tau \Vert_{\infty} <1$, i.e., system (\ref{psi_0}) has no solution. 
\end{proof}

We used the sparseness of the potential $B(x)$ in the Theorem (\ref{thm:sum B^beta inf}). What happens if $\sum\limits_{x \in Z^d}|B(x)|^{d/2} = \infty$ or $\sum\limits_{x \in Z^d}|B(x)|^{d/{\alpha}} = \infty$ but $B(x)$ has good estimate from below? The answer of this question depends on the construction of the correlation $B(x), \ x \in Z$. We will give one special case in the following theorem. 

\begin{theorem}
If we have \textbf{local} Laplacian operator $\Delta f(t,x)=\sum\limits_{\substack{z:|z-x|=1 \\ z\in Z^d}} (f(x+z)-f(x))$ and define $\mathcal{H}f(t,v)=\varkappa\Delta f(t,v)+B(v)f(t,v)$.	If $d \geq 3 $, $B(x) \geq \frac{C_0}{1 +|x|^2}$ and $C_0$ is fixed and large enough , and $\sum |B(x)|^{d/2} = \infty$, then $N_0(B) = \infty$. If $B(x)\geq \frac{C_1}{1+|x|^{\alpha}}$, $\alpha<2$, and $\sum |B(x)|^{d/\alpha} = \infty$, then $N_0(B)=\infty$ for arbitrary small $C_1$. 
\end{theorem}
\begin{proof}
Here we will give the proof for the first statement in this theorem, the analogy proof can be done for the second statement. 
Consider the cubes $Q_m=\{ x= (x_1, x_2, \cdots, x_d):  2^m < x_1  < 2^{m+1}$, $ 0 < x_i  < 2^{m}, i = 2,3, \cdots, d\}$. The cube $Q_m$ has volume size $2^{md}$. Then for $x\in Q_m$, the potential $B(x)$ has such lower bound 
	\begin{equation*}
		B(x) \geq \frac{C_0}{1 +|x|^2} \geq \frac{C_0}{1 + 2^{2(m+1)} + 2^{2m}+ \cdots + 2^{2m} } \geq 
		\frac{C_0}{1+(d+3)2^{2m}} \geq	\frac{C_0}{(d+3)2^{2m}}.
	\end{equation*}
	Then consider the following function for $ x= (x_1, x_2, \cdots, x_d) \in Q_m$ , 
	\begin{equation*}
		\psi_m(x_1, x_2, \cdots, x_d) = sin(k(x_1-2^m))	sin(kx_2)\cdots sin(kx_d).
	\end{equation*}
	Assume $\psi_m = 0$ on boundary of $Q_m$ ($\partial{Q_m}$), i.e., $sin(k2^{m}) = 0$  and $k = \frac{\pi}{2^m}$. 
	\begin{align*}
		\varkappa \Delta \psi_m(x) &= \varkappa\Delta_{x_1}\sin(k(x_1 -2^m))\sin(kx_2)\cdots \sin(kx_d)\\
		&+\varkappa\Delta_{x_2}\sin(k(x_1-2^m))\sin(kx_2)\cdots\sin(kx_d)+\cdots\\
		&+\varkappa\Delta_{x_d}\sin(k(x_1-2^m)\sin(kx_2)\cdots\sin(kx_d).
	\end{align*}
    
    Since $\varkappa\Delta\sin(\alpha x)=2\varkappa(\cos(\alpha)-1))\sin (\alpha x)$ for arbitrary $\alpha$,
	
	$$\lambda_m = 2 \varkappa(cos(k)-1) \sim -\frac{\varkappa \pi^2}{2^{2m}}.$$
	
	Then 
	$$\mathcal{L}\psi_m(x) \geq \lambda_m \psi_m +\frac{C_0}{(d
		+3)2^{2m}}\psi_m \sim [-\frac{\varkappa \pi^2}{2^{2m}} +\frac{C_0}{(d
		+3)2^{2m}}]\psi_m .$$
 
		Let the left hand side of preceding equation be positive, say $C_0 >  2\varkappa (d+3) \pi^2$,  then we will have  $$(\mathcal{L}\psi_m(x), \psi_m )\geq 0 $$ and we have infinitely many such cubes $Q_m$ and compactly supported test functions $\psi_m$ on $Q_m$, which means there are infinitely many positive eigenvalues. It is equivalently that if potential $B(x)$ decreases quadratically, there exist infinity many positive eigenvalues. 
\end{proof}

Then when we have infinitely many positive eigenvalues, the second moment $m_2(t,v)$ will diverge as $t\rightarrow \infty$. There is no steady state of $m_2$ in equation (\ref{m21_initial}).

Now we have answered the question about discrete positive spectrum for small $\varkappa$ if we have nonnegative $B$. We know that $\lambda_0(\mathcal{H}_2) >0$ for any $\varkappa > 0 $ if $B(x) \geq 0$ in dimension $d=1, 2$.  In the transient case, the positive eigenvalues are absent for large $\varkappa$, $\sum\limits_{z \in Z^d}|z|^2a(z) < \infty$ and $\sum\limits_{x \in Z^d}|B(x)|^{d/2} < \infty$ or in the case of the limiting stable laws with parameter $\alpha$,   $\sum\limits_{x \in Z^d}|B(x)|^{d/{\alpha}} < \infty$ (see Theorem  \ref{th: thereom 1}, \ref{th: thereom 2} ). We know from section 3 that $B(x)=\int_{T^d}e^{-ikx}\rho(k)dk$. It gives $B(0)=\int_{T^d}\rho(k)dk>0$.  $B(x)$ is a positively definite function  and it can be negative in some points and even it is possible that $B(x)<0,x\neq 0$. 

\textbf{Example 1}
Let $d=1$,and $\rho(k)=1-\sum\limits_{N=1}^{\infty}\frac{\cos(Nk)}{2^N}$, that is $B(x)=\begin{cases}
1\,\,&x=0\\
\frac{-1}{2^{1+|x|}}\,\,&x\neq 0
\end{cases}
$. But $\rho(k)\geq 0$, i.e., $B(x)$ is positively definite. Note that $\rho(0)=0$.

If the potential $B$ has negative and positive values, the positive eigenvalue $\lambda_0(\mathcal{H})$ can vanish due to the so called screening effect. The negative part of the the potential can overcome the attraction by the positive part.  Consider the following important example.

\textbf{Example 2}
Let $H\psi=\Delta \psi+\sigma V_0(x)\psi$, $\Delta \psi(x)=\psi(x+1)+\psi(x-1)-2\psi(x)$, $\sigma>0$ and 
$V_0(x)=\begin{cases}
1,\,\,x=0\\
-\frac{b}{2}\,\,x=\pm 1\\
0,\,\,|x|>1
\end{cases}$. Let us try to find the positive eigenvalue $\lambda_0(H)>0$ for the operator $H$ for small $\sigma$ and corresponding positive eigenfunction $\psi_0(x)>0$. We are looking $\psi_0(x)$ in the form with two unknown parameters $\mu,h>0$. 
\begin{align*}
\psi_0(0)=h,\psi_0(\pm1)=1,\psi_0(x)=\begin{cases}
e^{-\mu(x-1)},\,\,x\geq 1\\
e^{\mu(x+1)},\,\,x\leq -1
\end{cases}.
\end{align*}

Then solving the equation $H\psi_0=\lambda_0\psi_0$, we get the following results:
\begin{enumerate}[a).]
\item If $x\geq 2$ or $x\leq -2$, then 
\begin{align}\label{eq:|x|>1}
\lambda_0=2\cosh \mu-2.
\end{align} If $\sigma<<1$, then $\lambda_0<<1, \mu<<1$ and $\lambda_0\sim \mu^2$. 

\item If $x=0$, then 

\begin{align}\label{eq:|x|=0}
h=\frac{2}{\lambda_0-\sigma+2}.
\end{align}

\item If $x=1$ or $x=-1$, then 
\begin{align}\label{eq:|x|=1}
\lambda_0=h+e^{-\mu}-2-\frac{\sigma b}{2}.
\end{align}

It follows from (\ref{eq:|x|>1}),(\ref{eq:|x|=0}) and (\ref{eq:|x|=1}) that $h=\displaystyle e^{\mu}+\frac{\sigma b}{2}=\frac{1}{1+\frac{\lambda_0-\sigma}{2}}$. 
That is, when $b=1$,
\begin{align*} 
	1+(\sigma -\lambda_0/2) + (\sigma-\lambda_0/2)^2 + \cdots = e^{\mu}+\frac{\sigma b}{2}\\
		\frac{\sigma^2}{4} + \underline{\underline{o}}(\lambda_0 \sigma) = \mu + \underline{\underline{o}}(\mu^2).
\end{align*} 
Then
$\mu = \frac{\sigma^2}{4}+\underline{\underline{o}}( \sigma^2), \ \lambda_0= \mu^2 + \underline{\underline{o}}(\mu^2) = \frac{\sigma^4}{16} + \underline{\underline{o}}(\sigma^4) $.

\end{enumerate}

Our analysis indicates that the screening effect appears only if $b>1$. In the last case, $B(x)$ is not positively definite function. But if $b\leq 1$, that is $\sum\limits_{x\in Z^1}B(x)\geq 0$, then $\lambda_0(H)\sim C_1\sigma^2$ for $b>1$ and  $\lambda_0(H)\sim C_2\sigma^4$ for $b=1$, where $C_1,C_2$ are some constants and $B(x)$ is positively definite with $\rho(k)\geq 0$.

Let us stress that in the borderline situation $\sum\limits_{x\in Z^d} B(x)=0$ in the Example 2, the assumption of $\lambda_0(H)$ is different from the case when $\sum\limits_{x\in Z^d} B(x)>0$. The case when $\sum\limits_{x\in Z^d}B(x)>0$ has the positive eigenvalue for general potentials.

\begin{proposition}
Consider the spectral problem 
\begin{align*}
\mathcal{H}_2\psi(x) = 2\varkappa \Delta\psi +B(x) \psi = \lambda \psi(x),  
 x \in Z^1,  B(x) \in L^1(Z^1),  \sum_{z \in Z^1} B(x) = \delta> 0,
\end{align*}
then $\forall(\varkappa>0), \ \exists (\lambda_0(\mathcal{H}_2)) >0$.
\end{proposition}

\begin{proof}
Let's define the following even function $\psi_1(x) =\psi_1(-x) $
\begin{equation*}
\psi_1(x) = 
\begin{cases}
	1,  \ x \in [-L, L]\\
	\frac{L_1-x}{L_1-L}, \  x  \in [L, L_1]\\
	0, \ \text{otherwise}.
\end{cases}
\end{equation*}
Here $0 < L<L_1$ are two large parameters: $L \gg 1, \  L_1-L \gg 1$.

Let's Calculated the Dirichlet quadratic functional:
$
(\mathcal{H}_2 \psi_1, \psi_1) =2\varkappa(\Delta \psi_1, \psi_1) +(B(x) \psi_1, \psi_1)
$
Note that 
$$
(B(x) \psi_1, \psi_1)= \sum\limits_{x \in [-L, L]}B(x) +R_{L,L_1},  |R_{L,L_1}| \leq \sum_{|x|>L} |V(x)|.
$$
For any $0 < \delta_1 <\delta$, say $\delta_1 = \displaystyle\frac{1}{2} \delta$, one can find large enough $L$ such that $(B(x) \psi_1, \psi_1) \geq \delta_1$. But $\Delta \psi_1(x)=0$ if $x \neq \pm L$,  $x \neq \pm L_1$. And $|\Delta \psi_1(L)| =|\Delta \psi_1(L_1)| =\frac{1}{L_1-L} $. Finally,
\begin{equation*}
(\mathcal{H} \psi_1, \psi_1) \geq 1 \frac{4 \varkappa}{L_1-L} + \frac{1}{2} \sum_{x \in Z^d} V(x) = \frac{4 \varkappa}{L_1-L} + \frac{1}{2}  \delta \geq \frac{1}{3} \delta
\end{equation*}
if only for given $\varkappa$ and  $L$ the second parameter $L_1$ is the large enough. 
\end{proof}

This elementary proposition can be extended in several directions but it will be the subject of the different publications. 

\section*{Acknowledgments}
Dan Han was supported by University of Louisville EVPRI Grant "Spatial Population Dynamics with Disease" and AMS Research Communities "Survival Dynamics for Contact Process with Quarantine".

\section*{Declarations}
All authors certify that they have no affiliations with or involvement in any organization or entity with any financial interest or non-financial interest in the subject matter or materials discussed in this manuscript.

\bibliographystyle{plain}

\begin{thebibliography}{1}
\bibitem{carmona1994parabolic}
Ren{\'e}~A Carmona and Stanislav~A Molchanov, \emph{Parabolic anderson problem
  and intermittency}, vol. 518, American Mathematical Soc., 1994.
  
\bibitem{agbor2015global}
A~Agbor, S~Molchanov, and B~Vainberg, \emph{Global limit theorems on the
  convergence of multidimensional random walks to stable processes},
  Stochastics and Dynamics \textbf{15} (2015), no.~03, 1550024.


\bibitem{gartner1990parabolic}
J{\"u}rgen G{\"a}rtner and Stanislav~A Molchanov, \emph{Parabolic problems for
  the anderson model}, Communications in mathematical physics \textbf{132}
  (1990), no.~3, 613--655.

\bibitem{getan2017intermittency}
A~Getan, S~Molchanov, and B~Vainberg, \emph{Intermittency for branching walks
  with heavy tails}, Stochastics and Dynamics \textbf{17} (2017), no.~06,
  1750044.

\bibitem{molchanov2007quasicumulants}
 S.A. Molchanov, A.I. Petrov and N. Squartini, \emph{Quasicumulants and limit
  theorems in case of the cauchy limiting law}, Markov Process. Related Fields
  \textbf{13} (2007), no.~3, 597--624.

\bibitem{molchanov2010general}
Stanislav Molchanov and Boris Vainberg, \emph{On general
  cwikel--lieb--rozenblum and lieb--thirring inequalities}, Springer, 2010.

\bibitem{molchanov2019population}
Stanislav Molchanov and Boris Vainberg,\emph{Population dynamics with moderate tails of the underlying random
  walk}, SIAM Journal on Mathematical Analysis \textbf{51} (2019), no.~3,
  1824--1835.

\bibitem{Yarovaya}
Stanislav~Alekseevich Molchanov and Elena~Borisovna Yarovaya, \emph{Large
  deviations for a symmetric branching random walk on a multidimensional
  lattice}, Proceedings of the Steklov Institute of Mathematics \textbf{282}
  (2013), no.~1, 186--201.

\bibitem{zel1987intermittency}
Ya.~B. Zel'Dovich, S.A. Molchanov, A.~A. Ruzmaikin, and D.~D. Sokolov,
  \emph{Intermittency in random media}, Cambridge Scientific Publishers
  Limited, 2014.

\end{thebibliography}

\end{document}